\documentclass[oneside,a4paper,11pt]{article}
\usepackage{amsmath}
\usepackage{amssymb}
\usepackage{graphicx}
\usepackage{natbib}
\newtheorem{Definition}{\bf \large Definition}[section]
\newtheorem{Theorem}{\bf \large Theorem}[section]
\newtheorem{Lemma}[Theorem]{\bf \large Lemma}

\newtheorem{Remark}{\bf \large Remark}[section]
\renewcommand{\theequation}{\thesection.\arabic{equation}}
\newcommand{\proof}{\noindent \textbf{Proof.} \hspace {0.3cm}}

\newcommand{\pr}{\mathrm{P}}
\newcommand{\E}{\mathrm{E}}

\newcommand{\dif}{\mathrm{d}}

\begin{document}

\title{Optimal dividend  control for a generalized risk model with investment incomes and debit interest}
\author{Jinxia Zhu \\Actuarial Studies\\
 Australian School of Business\\ University of New South Wales\\
  Australia}
\date{ }

\maketitle
\makeatletter{\renewcommand*{\@makefnmark}{}
\footnotetext{Actuarial Studies, Australian School of Business, University of New South Wales, NSW 2052, Australia}\makeatother}
\makeatletter{\renewcommand*{\@makefnmark}{}
\footnotetext{Email: jinxia.zhu@unsw.edu.au}\makeatother}

\begin{abstract}
This paper investigates dividend optimization of an insurance corporation under a more realistic model which takes into consideration refinancing or capital injections. The model follows the compound Poisson framework with credit interest for positive reserve, and debit interest for negative reserve. Ruin occurs when the reserve drops below the critical value.  The company controls the dividend pay-out dynamically with the objective to maximize the expected total discounted dividends until ruin. We show that that the optimal strategy is a band strategy and it is optimal to pay no dividends when the reserve is negative.
\end{abstract}

 \noindent {\bf Key words}
  Absolute ruin,  dividend optimization,   stochastic control, value function, viscosity solution.

\bigskip

\noindent  {\bf Mathematics Subject Classification (2000)}  91B30; 93E20; 49L25\\
\noindent {\bf JEL Classification} C61; C02

\section{Introduction}
Dividend optimization problems for financial and insurance corporations have attracted extensive attention over the last few decades. One  of this type of  problems is to find the optimal dividend pay-out scheme, i.e. choosing the times and amounts of dividend payments to maximize the objective function - the expected total discounted dividend pay-outs until the time of ruin.

 In the area of non-life insurance, a well established model  for the cash reserve  is the Cram\'er-Lundberg model (also called the compound Poisson model or the classical risk model), which is based on Poisson claim-arrivals and linear premium income. \cite{EmbrechtsSchmidli1994} claimed that ``many of the `rules of thumb' used in practice can be traced back to the classical Cram\'er-Lundberg model". However, starting from the middle of 1990's,  a large number of papers dealing with optimization problems for insurance companies, use  the diffusion process - a limiting process of the Cram\'er-Lundberg model, to model the reserve in the absence of the dividends, e.g. \cite{Jeanblanc-PicqueShiryaev1995}, \cite{CadenillasChoulliTaksarZhang2006} and \cite{Paulsen2007}. Diffusion process modeling of the reserve process allows the use of optimal diffusion control techniques and is therefore more mathematically tractable. A survey of optimal dividend control for diffusion processes can be found in \cite{Taksar2000}.

  There have been a few attempts to study the dividend optimization problem under the Cram\'er-Lundberg model.  \cite{Gerber1969} considered the  dividends optimization problem for a  classical Cram\'er-Lundberg model   and  proved  that the corresponding optimal dividend strategy is a band strategy. \cite{AzcueMuler2005} considered the  Cr\'amer-Lundberg model with reinsurance and dividend payments and proved the optimal dividend payment policy maximizing the expected total discounted dividend pay-outs  is also a band strategy. \cite{AlbrecherThonhauser2008} studied the dividend optimization problem in the  Cr\'amer-Lundberg setting including constant force of interest and pointed out that the optimal strategy is also of band type. \cite{KulenkoSchmidli2008} found that the optimal dividend strategy for the   Cr\'amer-Lundberg model with capital injections is a barrier strategy. For applications of stochastic control in insurance, please refer to \cite{Schmidli2008} and references therein. A list of literature on dividend optimization problems under the Cram\'er-Lundberg model  can be found in \cite{AlbrecherThonhauser2009}. For a review of  dividend strategies in the actuarial literature, see  \cite{Avanzi2009}.

 There has been extensive work dedicated to the generalization of the classical risk model to suit more realistic situations. One way of generalization is to  allow the company to refinance when the company is in deficit and the deficit is not too large.  The idea was developed by  \cite{Borch1969}, where he proposed that ruin (negative reserve) does not mean the end of game but only the necessity of raising additional money. He argued that ``insurance companies get into difficulties fairly regularly and rescue operations are considered in the insurance world, if not daily, at least annually" and that it will be a good investment to rescue a company when the situation is not too serious, and concluded that a company should be rescued if the benefits exceed the cost of the new financing required, e.g. when the deficit is not too large.  Since then the ``absolute ruin model" has been developed, where the company is allowed to borrow money to settle the claims  if the reserve is negative but still above the critical level so that it can continue its business.   The company will need to pay interest (debit interest) on the loan  and  pay back  debt  interest continuously from the received premiums. The critical level is the value of reserve below which the premiums received are insufficient to cover interest payments on the debt. Absolute ruin occurs when the reserve reaches or drops below the critical level for the first time.

The absolute ruin problem has received considerable attention.
 \cite{Gerber1971} studied the absolute ruin probability in the compound Poisson model.  \cite{EmbrechtsSchmidli1994} considered the absolute ruin probability when the reserve process is a piecewise-deterministic Markov process. \cite{DicksonEgidiodosReis1997} used simulation to study the Cram\'er-Lundberg model with absolute ruin.  \cite{CaiGerberYang2006} studied an Ornstein-Uhlenbeck type model with credit and debit interest. \cite{Cai2007} discussed the Gerber-Shiu function in the classical risk model with absolute ruin. \cite{GerberYang2007} investigated the absolute ruin probability based on the classical risk model perturbed by diffusion with investment. \cite{ZhuYang20082} studied the asymptotic behavior of the absolute ruin probability in the Cram\'er-Lundberg model with credit and debit interest. Some other related references are \cite{YuenZhouGuo2008} and \cite{WangYin2009}.

In this paper, we consider the dividend optimization under the the compound Poisson model with credit interest for positive reserve, and debit interest for negative reserve. The paper is organized as follows. Section 2 presents the model and formulates the dividend optimization problem. In section 3,  we derive some basic and important
properties of the  value function, and characterize  the value function as the unique nonnegative and nondecreasing viscosity super-solution of  the associated
Hamilton-Jacobi-Bellman equation  that satisfies a linear growth condition and a boundary condition. In section 4, we prove the existence of the optimal dividend strategy and identify the optimal dividend payout scheme as a band strategy. It is shown that the optimal strategy is to pay no dividends when the reserve is negative. A conclusion is provided in section 5.
\section{The model and the optimization problem}
  Consider a continuous
time  model for the surplus of an insurance company where claims arrive according to a Poisson process with intensity rate $\lambda$ and premiums are
collected continuously at the rate $p$. The amount of each
claim is independent of its arrival time, and is also independent of any other claims.  Let $S_i$ denote
the arrival epoch of the $i$th claim
 and $U_i$ its size. Let $N(t)=\sharp\{i:S_i\leq t\}$. Then $N(t)$ is the number of claims up to time
 $t$ and  follows a Poisson process with  rate $\lambda$. The sequence $\{U_i\}$ is assumed to be identically
 and independently distributed with distribution function $F(\cdot)$ and   independent of
 $\{N(t)\}$. Moreover, the insurance company earns credit
interest under a constant force $r$ ($r>0$) when the surplus is
positive, and when the surplus drops below $0$, the insurer could borrow money
with the amount equal to the deficit  under force of debit interest
  $\alpha>r$. In the mean time, the insurer will repay the
  debts and the debt interest continuously from the premium incomes.
  This leads to the  following dynamics for the risk reserve process $\{X_t\}_{t\geq 0}$ in the absence of dividend payments: $$\dif X_t=(p+r X_{t-}I\{X_{t-}\geq 0\}+\alpha X_{t-}I\{X_{t-}<0\})\dif t
-\dif Y_t,$$
where $X_t$ represents the surplus at time $t$ and $ Y_t=\sum_{i=1}^{N(t)}U_i$ is the aggregate claim up to
time $t$.

Now suppose the company pays dividends to its shareholders with the
accumulative amount of dividends paid up to time $t$ being denoted by  $L_t$. Let $\overline{R}^L_t$  denote the controlled reserve at time $t$. Then
\begin{eqnarray}\dif \overline{R}^L_t=(p+r \overline{R}^L_{t-}I\{\overline{R}^L_{t-}\geq 0\}+\alpha \overline{R}^L_{t-}I\{\overline{R}^L_{t-}<0\})\dif t
-\dif \left(\sum_{i=1}^{N(t)}U_i\right)-\dif L_t.\label{dyn}
\end{eqnarray}
The company controls dynamically the dividend pay-outs: the times and the amounts of dividends to be paid out. A control
strategy is described by a dividend distribution process $L=\{L_t\}_{t\geq 0}$.

 Notice from
the above dynamics that the premium incomes will no longer be able to cover
the debts when the surplus is less than or equal to $-\frac p\alpha $. That is,
the surplus process will not be able to return to a positive amount
whenever the process hits $-\frac p\alpha $ or any level below that.
We call $-\frac p\alpha $ the critical value and  define
\textit{the time of ruin} as $T^L=\inf\{t\geq 0: \overline{R}^L_t\leq -\frac p\alpha\}.$
  The  time of ruin defined above is also called the  time of absolute ruin in the sense that the
  surplus will no longer be able to return to a positive level.

  All our random quantities are defined on the   complete probability space $(\Omega,\mathcal{F},\pr)$. Let $\mathcal{N}$ denote the class of null sets in $\Omega$ and define $\mathcal{F}_t=\sigma(X_0,Y_s,0\leq s\leq t)\bigvee \mathcal{N}$. Throughout the paper, we base our study on the filtered probability space  $(\Omega,\mathcal{F},\{\mathcal{F}_t\}_{t\geq 0},\pr)$.

A control strategy is
\textit{admissible} if the
process    $\{L_t\}_{t\geq 0}$  with $L_0=0$, is predictable, nondecreasing, left
continuous with right limits (c\'agl\'ad) and satisfies 
the requirement  that paying dividends would not cause ruin immediately.
We use $\Pi$ to denote the set of all admissible strategies.

Define $E_x[\ \cdot \ ]=\E[\ \cdot\ |\overline{R}_0=x].$
Let $\delta$ be the force of discount with $\delta>r$.
Given the initial reserve $x$, the performance of a dividend strategy $L$ is measured by  the expectation of the
cumulative discounted dividends until ruin, i.e.
\begin{eqnarray}
V_L(x)=\E_x\left[\int_0^{T^L} e^{-\delta s}\dif L_s\right]. \label{7a1}
\end{eqnarray}
The integral here is interpreted path-wise in a Lebesgue-Stieltjes
sense. The function $V_L(x)$ is  called the \textit{return function}. Obviously, $V_L(x)=0\ \ \mbox{ for $x\leq -\frac p\alpha$}.$

 The objective of the company
is to find an optimal dividend payout scheme $L$ in the set of admissible strategies $\Pi$ such that the
expectation of total discounted dividend pay-outs until the time of ruin is maximized.

Define the \textit{value function} (also called the \textit{ optimal return function}) by $V(x)=\sup_{L\in \Pi}V_L(x).$
 If there exists a control strategy $L^\ast$ such that $V(x)=V_{L^\ast}(x)$, then $L^\ast$ is called the optimal dividend distribution process (the optimal dividend strategy).

It can be seen that $T^L$ is a stopping
 time.
 In the paper,  we will consider the stopped process
$R^L_t=\overline{R}^L_tI\{t<T^L\}-\frac p\alpha I\{t\geq T^L\}.$

 To simplify the notation we will omit the superscripts $L$ in $T^L$ and $R^L$.

Since the reserve process  in the absence of the control variable is a Markov process, the problem here is the optimization problem for  a controlled Markov process.
As the cumulative dividend process $L$  may not be continuous with respect to time, the optimization problem  is a singular control problem. In the context of stochastic control theory, the optimization problem can be associated with a Hamilton-Jacobi-Bellman (HJB) equation derived by using the Dynamic Programming Principle. In this case, the HJB equation is a first-order integro-differential equation. However, the differentiability of the value function is a question. Actually, even under a specifically predetermined dividend strategy, the differentiability of the corresponding return function can not be guaranteed. It was shown in \cite{ZhuYang2009} that the differentiability of the return function under a barrier or threshold dividend strategy depends on the level of smoothness of the claim size distributions.  In this paper, we show that  the value function is absolutely continuous but may not be differentiable. So we resort to the concept of viscosity solutions.

Based on  techniques of probability and Stochastic Control theory, we show that the value function is a viscosity solution of the associated HJB equation and it is the unique solution  satisfying certain regularity and boundary conditions. We also prove that the optimal dividend payment strategy exists and is of a band type, an and that it is optimal to pay no dividends at all when the surplus is negative.
Proofs of some lemmas and theorems are relegated to the appendix.
\section{The value function}
In this section, we derive some analytical  properties of the value
function $V(x)$. We show that $V(x)$ is not necessarily differentiable everywhere, but almost everywhere, and that the value function is the viscosity solution to the associated HJB equation but not necessarily the classical solution. It will also be proven that the value function is the unique solution satisfying certain conditions.
\begin{Theorem}\label{thm2.1} If $r<\delta$, $V(x) \geq x+\frac p\alpha$  for $x\in\mathbb{R}$, and  $V(x)\leq
\frac{\delta x+p}{\delta-r}+\frac p\alpha$ for $x\geq 0$.
 \end{Theorem}
\proof
 To prove the lower bound, consider a dividend payout scheme such that the part of initial reserve in excess of the critical value $-\frac p\alpha$ is paid out immediately as dividends.  Then ruin occurs immediately. In this case,  the return function given the initial reserve $x$, is $x+\frac p\alpha$. So the optimal return function $V(x)$ is always greater than or equal to $x+\frac p\alpha$.

From \eqref{dyn} we can see that given that the initial reserve is nonnegative,  the  inequality $\dif R_{t}\leq (p+r R_{t-})\dif t$ holds. As a result, given $R_0=x$ we have
$
R_t\leq
 \frac 1r\left(e^{rt}(p+rx)-p\right)$ for $ x\geq 0$. Hence,
by integration by parts,  $L_t\leq R_t+\frac p\alpha$ for any $L\in\Pi$ and the definition \eqref{7a1}  we can obtain
 $
V_L(x)=\E_x\left[\int_0^T e^{-\delta s}\dif L_s\right]\leq \E_x\left[\int_0^\infty \delta L_s e^{-\delta s}\dif s\right]\label{dd1}=\frac{\delta x+p}{\delta-r}+\frac p\alpha$ for  $x\geq 0$.
 \hfill $\square$\\

 Define
\begin{eqnarray}
t_0(x,y)=\begin{cases} \frac 1r \log(\frac{ry+p}{rx+p}) & y>x\geq 0\\
                       \frac 1r \log(\frac{r y+p}{p})+ \frac 1\alpha \log(\frac{p}{\alpha x+p})& y> 0>x>-\frac p\alpha\\
                       \frac 1\alpha \log(\frac{\alpha y+p}{\alpha x+p})& 0\geq y>x>-\frac p\alpha
         \end{cases}.\label{cc2}
\end{eqnarray}
The quantity $t_0(x,y)$ is equivalent to the time it takes for the surplus process with initial value $x$ to reach $y$ $(y>x)$ for the first time given that there are no claims and no dividends paid out.

\begin{Theorem}\label{thm2.2}The value function $V$ satisfies the following inequalities \begin{eqnarray}
y-x\leq V(y)-V(x)\leq V(x)\begin{cases}\left((\frac{ry+p}{rx+p})^{\frac{\lambda+\delta}{r}}-1\right)& y>x\geq 0\nonumber\\
                                       \left((\frac{r y+p}{p})^{\frac{\lambda+\delta}{r}}(\frac{p}{\alpha x+p})^{\frac{\lambda+\delta}{\alpha}}-1\right)& y\geq 0>x>-\frac p\alpha\nonumber\\
                                       \left((\frac{\alpha y+p}{\alpha x+p})^{\frac{\lambda+\delta}{\alpha}}-1\right)& 0>y>x>-\frac p\alpha
                                       \end{cases}.\label{e02}
\end{eqnarray}
\end{Theorem}
\proof (i) We first prove the lower bound. For any $\epsilon>0$, let
$L_\epsilon(x)$ denote an admissible $\epsilon$-optimal strategy given the initial
reserve $x$, i.e. $V_{L{\epsilon(x)}}(x)\geq V(x)-\epsilon.$

For $y>x>-\frac p\alpha$, given the initial reserve $R_0=y$ we use $L(y,x)$ to denote a strategy that pays
an amount $y-x$ as dividends immediately and then pays dividends according to  the strategy
$L_\epsilon(x)$. Then given the initial reserve $R_0=y>x$, under the strategy $L(y,x)$ we have $V_{L(y,x)}(y)=y-x+V_{L_\epsilon(x)}(x).$
So for any $\epsilon>0$,
$
V(y)\geq y-x+V_{L_\epsilon(x)}(x)\geq y-x+V(x)-\epsilon.$
Consequently, $V(y)-V(x)\geq y-x.$

(ii) Now, we proceed to prove the upper bounds.  For $y>x>-\frac
p\alpha$, for the surplus process with initial reserve $x$, let $\tau(x,y)$ denote the time it will take for the surplus process to reach up to $y$ for the first time, and  define the strategy
$\hat{L}(x,y)$ as follows:
\begin{itemize}
\item   pay out no dividends until the
reserve reaches $y$,
 \item  then at the moment that the reserve reaches $y$ for the first time ($\tau(x,y)$), treat the reserve process as a new
process that starts at this moment with initial capital $y$, and apply
the strategy
$L_\epsilon(y)$, i.e. $\theta_{\tau(x,y)}\hat{L}(x,y)=L_\epsilon(y)$.

\end{itemize}

 Note that starting from the initial value $x>-\frac p\alpha$, ruin will not occur before the arrival of the first claim ($S_1$), and the reserve will
 reach $y$ ($y>x$) at time $t_0(x,y)$ if no claims arrive before time $t_0(x,y)$, that is $\tau(x,y)=t_0(x,y) \ \ \mbox{ on } \{S_1>t_0(x,y)\}.$
Then for $y>x>-\frac p\alpha$ and for $\epsilon>0$, by noticing that $V_{L_\epsilon(y)}(y)\geq V(y)-\epsilon$ we have
\begin{eqnarray*}
V(x)&\geq& V_{\hat{L}(x,y)}(x)= \E_x[e^{-\delta \tau(x,y)}V_{L_\epsilon(y)}(y);\tau(x,y)\leq T]\\
&\geq& \E_x[e^{-\delta \tau(x,y)}V_{L_\epsilon(y)}(y);S_1>t_0(x,y)]\geq  e^{-(\lambda+\delta) t_0(x,y)}(V(y)-\epsilon).
\end{eqnarray*}
Hence, $V(y)-V(x)\leq V(x)(e^{(\lambda+\delta) t_0(x,y)}-1).$
This combined with \eqref{cc2} gives the upper bounds.
 \hfill $\square$\\

\begin{Theorem}\label{thm2.3} The value function $V(x)$ is nonnegative, nondecreasing, continuous on $[-\frac p\alpha,
  \infty)$ and locally Lipschitz continuous on $(-\frac p\alpha,\infty)$. Therefore, $V^\prime(x)$ exists almost
  everywhere on  $(-\frac p\alpha,\infty)$. Furthermore,  $V^\prime(x)\geq 1$, if $V^\prime(x)$ exists.\end{Theorem}
  \proof
  All the stated properties of $V(x)$ are direct results of Theorem \ref{thm2.1} and Theorem \ref{thm2.2} except for  the right continuity of $V(x)$ at $x=-\frac p\alpha$.

  To prove the right continuity, it is sufficient to show that $\limsup_{x\downarrow -\frac p\alpha}V(x)=0$. If this is not true, then we can find a sequence $\{x_n\}$ with $x_n\downarrow -\frac p\alpha$ such that $\lim_{n\rightarrow\infty}V(x_n)>0$, that is, there exists an $\epsilon_0>0$ and $N$ such that $  V(x_n)>\epsilon_0$ for all $n\geq N$.
  Let $L^{(x,\frac{\epsilon_0}2)}$ denote a $\frac{\epsilon_0} 2$-optimal strategy for the reserve process with  initial reserve  $x$, that is, $V_{L^{(x,\frac{\epsilon_0}2)}}(x)\geq V(x)-\frac{\epsilon_0}2.$
Then, we have
\begin{eqnarray}
V_{L^{(x_n,\frac{\epsilon_0}2)}}(x_n)\geq V(x_n)-\frac{\epsilon_0}2>\frac{\epsilon_0}2 \ \ \mbox{ for $n\geq N$}.\label{cb03}
\end{eqnarray}

Consider a stochastic process $\{R_t^\prime\}$ with dynamics
$\dif R_t^\prime= (p+r R_{t}^\prime I\{R_{t}^\prime\geq 0\}+\alpha R_{t}^\prime I\{R_{t}^\prime<0\})\dif t.$
Given $R_0^\prime=x$, integration yields
\begin{eqnarray}
R_t^\prime\leq
 \begin{cases}
 \frac 1r\left(e^{rt}(p+rx)-p\right)&  x\geq 0\\
  \frac 1\alpha\left(e^{\alpha t}(p+\alpha x)-p\right)&\  x< 0, t\leq t_0(x,0)\\
  \frac 1r\left(e^{r(t-t_0(x,0))}p-p\right)&  x< 0,\ t> t_0(x,0)
 \end{cases}.\label{dd2}
\end{eqnarray}

Note that $R_t\leq R_t^\prime$ given that $R_0=R_0^\prime>-\frac p\alpha$. Using the fact that $L_t\leq R_t+\frac p\alpha\leq R_t^\prime+\frac p\alpha$ and \eqref{dd2}, by integration by parts it follows from \eqref{7a1} that   for $x\in(-\frac p\alpha,0)$
\begin{eqnarray}
V_L(x)
&\leq& \E_x[\int_0^\infty \delta L_s e^{-\delta s}\dif s]\nonumber\\
&\leq& \delta \int_0^{t_0(x,0)} \frac1\alpha\left(
e^{\alpha s}(p+\alpha x)-p\right)e^{-\delta s}\dif s\nonumber\\&&+\delta\int_{t_0(x,0)}^\infty \frac 1r\left(e^{r(s-t_0(x,0))}p-p\right)e^{-\delta s}\dif s+\frac p\alpha.\label{cd01}
\end{eqnarray}
Notice that the expression on the right-hand side of \eqref{cd01} has limit $0$ as $x\downarrow -\frac p\alpha$ and does not depend on $L$. So we can find an $N^\prime$ such that for all $n\geq N^\prime$, $V_L(x_n)<\frac{\epsilon_0}4$ holds  for all admissible strategy $L$.
Therefore, setting $L$ to be $ L^{(x_n,\frac{\epsilon_0}2)}$ gives $V_{L^{(x_n,\frac{\epsilon_0}2)}}(x_n)<\frac{\epsilon_0}4\ \ \mbox{ for all $n\geq N^\prime$},$
which is a contradiction to \eqref{cb03}.
Hence, the value function $V(x)$ is right continuous at $-\frac p\alpha$.
  \hfill $\square$\\

\bigskip

Applying standard arguments from stochastic control  theory (e.g. \cite{FlemingSoner1993}) or an approach
similar  to that in \cite{AzcueMuler2005}, we can show that the optimal value
function fulfils the Dynamic Programming Principle:
 $$V(x)=\sup_{L\in\Pi}\E_x\left[\int_0^{\tau\wedge T}e^{-\delta s}\dif
 L_s+e^{-\delta (\tau\wedge T)}V(R_{\tau\wedge T})\right
]
 \ \ \mbox{ for any
   stopping time $\tau$,}$$
and the associated Hamilton-Jacobi-Bellman (HJB) equation is
\begin{eqnarray}
\max\{1-V^\prime (x),\mathcal{L}_V(x)\}=0,\label{1}
\end{eqnarray}
where $\mathcal{L}$ is a generator defined by \begin{eqnarray}
\mathcal{L}_V(x)&=&\left(p+rxI\{x\geq 0\}+\alpha x I\{x< 0\}\right) V^\prime(x)\nonumber\\
&&-(\lambda+\delta)V(x)+\lambda\int_0^{x+\frac p\alpha }V(x-u)\dif
F(u).\label{3}
\end{eqnarray}

Although from the last section we know that  $V^\prime(x)$ exists almost everywhere, we have no guarantee that $V(x)$ is differentiable for all $x>-\frac
p\alpha$. Therefore, we can not  expect $V(x)$ to be a classical solution
to the HJB equation. In the following we will show that the
value function $V(x)$ is a  viscosity solution to the HJB equation \eqref{1}, and that $V(x)$ is the unique nonnegative,  nondecreasing and locally Lipschitz continuous viscosity solution of \eqref{1} satisfying a linear growth condition and the  boundary condition $V(-\frac p\alpha)=0$.

\begin{Definition}
(i) A continuous function $\underline{u}:[-\frac
p\alpha,\infty)\rightarrow \mathbb{R}$ is said to be a viscosity
sub-solution   of \eqref{1} on $(-\frac
p\alpha,\infty)$ if for any $x\in (-\frac
p\alpha,\infty)$ each continuously differentiable function
$\phi:(-\frac p\alpha,\infty)\rightarrow \mathbb{R}$ with
$\phi(x)=\underline{u}(x)$ such that $\underline{u}-\phi$ reaches the
maximum at $x$ satisfies
$$\max\{1-\phi^\prime(x),\mathcal{L}_\phi (x)\}\geq 0.$$\\
\noindent (ii) A continuous function $\overline{u}:[-\frac
p\alpha,\infty)\rightarrow \mathbb{R}$ is said to be a viscosity
super-solution   of \eqref{1} on $ (-\frac
p\alpha,\infty)$ if for any $ x\in(-\frac
p\alpha,\infty)$ each  continuously differentiable function
$\phi:(-\frac p\alpha,\infty)\rightarrow \mathbb{R}$ with
$\phi(x)=\overline{u}(x)$ such that $\overline{u}-\phi$ reaches the
minimum at $x$ satisfies
$$\max\{1-\phi^\prime(x),\mathcal{L}_\phi (x)\}\leq 0.$$
\noindent (iii)  A continuous function ${u}:[-\frac
p\alpha,\infty)\rightarrow \mathbb{R}$ is  a viscosity
solution   of \eqref{1} on $ (-\frac
p\alpha,\infty)$ if it is both a viscosity sub-solution and a
viscosity super-solution on $(-\frac
p\alpha,\infty)$. \label{Def}
\end{Definition}

 For any continuously differentiable function $\phi$ and any continuous function $v$, define an operator
$\mathcal{L}_{v,\phi}(x)=\left(p+rxI\{x\geq 0\}+\alpha x I\{x< 0\}\right) \phi^\prime(x)-(\lambda+\delta)v(x)+\lambda\int_0^{x+\frac p\alpha }v(x-u)\dif
F(u).$ As has been shown in  (\cite{Sayah1991} and \cite{BenthKarlsenReikvam2001}), the definition of viscosity sub and super solutions has the following alternative version.
\begin{Definition}\label{alt}
(i) A continuous function $\underline{u}:[-\frac
p\alpha,\infty)\rightarrow \mathbb{R}$ is said to be a viscosity
sub-solution   of \eqref{1} on $(-\frac
p\alpha,\infty)$ if for any $x\in(-\frac
p\alpha,\infty)$ each continuously differentiable function
$\phi:(-\frac p\alpha,\infty)\rightarrow \mathbb{R}$ with
$\phi(x)=\underline{u}(x)$ such that $\underline{u}-\phi$ reaches the
maximum at $x$ satisfies
\begin{eqnarray*}\max\{1-\phi^\prime(x),\mathcal{L}_{\underline{u},\phi}(x)\}\geq 0.\end{eqnarray*}
\noindent (ii) A continuous function $\overline{u}:[-\frac
p\alpha,\infty)\rightarrow \mathbb{R}$ is said to be a viscosity
super-solution   of \eqref{1} on $(-\frac
p\alpha,\infty)$ if for any $(-\frac
p\alpha,\infty)$ each continuously differentiable function
$\phi:(-\frac p\alpha,\infty)\rightarrow \mathbb{R}$ with
$\phi(x)=\overline{u}(x)$ such that $\overline{u}-\phi$ reaches the
minimum at $x$ satisfies
$$\max\{1-\phi^\prime(x),\mathcal{L}_{\overline{u},\phi}(x)\}\leq 0.$$
\end{Definition}

The following remarks are standard in the context of viscosity theory (eg \cite{Capuzzo-DolcettaLions1990} and \cite{CrandallEvansLions1984}),
which will be useful in the
proof of our main results.
\begin{Remark}\label{cc50}
(i) For any viscosity sub-solution $\underline{u}$ on $(-\frac p\alpha,\infty)$, there exists a continuously differentiable function
$\phi:(-\frac p\alpha,\infty)\rightarrow \mathbb{R}$ such that $\underline{u}-\phi$
reaches a maximum at $x>-\frac p\alpha$ with $\phi^\prime(x)=q$ if and only if
$$\liminf_{y\uparrow
x}\frac{\underline{u}(y)-\underline{u}(x)}{y-x}\geq q\geq \limsup_{y\downarrow
x}\frac{\underline{u}(y)-\underline{u}(x)}{y-x}.$$
(ii)  For any viscosity super-solution $\overline{u}$ on $(-\frac
p\alpha,\infty)$, there exists a continuously differentiable function
$\phi:(-\frac p\alpha,\infty)\rightarrow \mathbb{R}$ such that $\overline{u}-\phi$
reaches a minimum at $x>-\frac p\alpha$ with $\phi^\prime(x)=q$ if and only if
$$\liminf_{y\downarrow
x}\frac{\overline{u}(y)-\overline{u}(x)}{y-x}\geq q\geq \limsup_{y\uparrow
x}\frac{\overline{u}(y)-\overline{u}(x)}{y-x}.$$
\end{Remark}

For any $t\geq 0$, define a functional $\mathcal{M}_t$ by
\begin{eqnarray}
\mathcal{M}_t(\phi)
&=&\sum_{s\leq t,R_{s-}\neq R_s}\left(\phi(R_s)-\phi(R_{s-})\right)e^{-\delta  s}\nonumber\\
&&-\lambda \int_0^te^{-\delta  s}\dif s\int_0^{\infty}\left(\phi(R_{s-}-y)-\phi(R_{s-})\right)\dif F(y).\label{cb04}
\end{eqnarray}
Then $\{\mathcal{M}_t(\phi)\}$ is a local martingale. If $\phi(\cdot)$ is bounded by a linear function, then $\mathcal{M}_t(\phi)$ is bounded below and therefore a super-martingale by applying Fatou's Lemma.

   Consider any  nonnegative and nondecreasing function $\phi$  and  any stopping time $\tau$ such that   ${\phi}^\prime(R_{t})$ exists for  all $t\leq \tau$ and \begin{eqnarray}\phi^\prime(R_t)\geq 1 \ \mbox{for all $t\leq \tau$}.\label{7a11}
   \end{eqnarray}

   Let  $\{L_t^c\}$ denote the continuous part of $\{L_t\}$. It can be seen that
\begin{eqnarray}
&& {\phi}(R_{\tau})e^{-\delta \tau}- {\phi}(R_0)=\int_0^{\tau}\dif \left( {\phi}(R_t)e^{-\delta t}\right)\nonumber\\
&=&\int_0^{\tau}  {\phi}^\prime(R_t)e^{-\delta t}\dif R_t-\delta\int_0^{\tau} {\phi}(R_t)e^{-\delta t}\dif t\nonumber\\
&=&\int_0^{\tau}  {\phi}^\prime(R_t)e^{-\delta t}\left(p+rR_tI\{R_t\geq 0\}+\alpha R_tI\{R_t< 0\}\right)\dif t\nonumber\\
&&-\int_0^{\tau}  {\phi}^\prime(R_t)e^{-\delta t}\dif L_t^c+\sum_{t\leq\tau,R_{t-}\neq R_t}( {\phi}(R_t)- {\phi}(R_{t-}))e^{-\delta t}\nonumber\\
&&+\sum_{t<\tau, R_t\neq R_{t+}}( {\phi}(R_{t+})- {\phi}(R_{t}))e^{-\delta t}-\delta\int_0^{\tau} {\phi}(R_t)e^{-\delta t}\dif t,\label{23}
\end{eqnarray}
where the last equality follows from the fact that $L_t$ is left-continuous and nondecreasing.

\noindent Since $R_t\neq R_{t+}$ only occurs at the jumps of $L_t$ and $L_t$ is left-continuous in $t$, then $ R_{t+}-R_t=-(L_{t+}-L_t)$
and
\begin{eqnarray}\sum_{t<\tau, R_t\neq R_{t+}}( {\phi}(R_{t+})- {\phi}(R_{t}))e^{-\delta t}=-\sum_{t<\tau, R_t\neq R_{t+}}e^{-\delta t}\int_0^{L_{t+}-L_t}  {\phi}^\prime(R_{t}-u)\dif u.\label{08}
\end{eqnarray}
  Then by \eqref{08} and \eqref{7a11} we have
\begin{eqnarray}&&-\int_0^{\tau} {\phi}^\prime(R_t)e^{-\delta t}\dif L_t^c+ \sum_{t<\tau, R_t\neq R_{t+}}( {\phi}(R_{t+})- {\phi}(R_{t}))e^{-\delta t}\nonumber\\
&\leq& -\int_0^{\tau}e^{-\delta t}\dif L_t^c-\sum_{t<\tau, R_t\neq R_{t+}} e^{-\delta t}\left(\int_0^{L_{t+}-L_t}\dif u\right)\nonumber\\
&=&-\int_0^{\tau}e^{-\delta t}\dif L_t.\label{24}
\end{eqnarray}
Using \eqref{3}, \eqref{cb04}, \eqref{23} and \eqref{24} and noting that $\phi(x)\geq 0$ for $x\leq -\frac p\alpha$, we have
\begin{eqnarray}
&& {\phi}(R_{\tau})e^{-\delta \tau}- {\phi}(R_0)\nonumber\\
&\leq&\int_0^{\tau}  {\phi}^\prime(R_{t-})e^{-\delta t}\left(p+rR_{t-}I\{R_{t-}\geq 0\}+\alpha R_{t-}I\{R_{t-}< 0\}\right)\dif t\nonumber\\
&&-\int_0^{\tau} e^{-\delta t}\dif L_t+\mathcal{M}_{\tau}( {\phi})+\lambda \int_0^{\tau}e^{-\delta t}\dif t\int_0^{\infty}\left( {\phi}(R_{t-}-u)- {\phi}(R_{t-})\right)\dif F(u)\nonumber\\
&&-\delta\int_0^{\tau} {\phi}(R_t)e^{-\delta t}\dif t\nonumber\\
&=&\int_0^{\tau}\mathcal{L}_ {\phi}(R_{t-})e^{-\delta t}\dif t-\int_0^{\tau} e^{-\delta t}\dif L_t+\mathcal{M}_{\tau}( {\phi}).\label{26}
\end{eqnarray}

\bigskip

In the next theorem, we  show that the value function $V$
is a viscosity solution of the HJB equation \eqref{1}.
\begin{Theorem}\label{thm3.1}
(i) $V(x)$ is a viscosity solution of  \eqref{1} on $(-\frac p\alpha,+\infty)$.\\
(ii)  Define $D_V(x,h_n) = \frac{V(x+h_n)-V(x)}{h_n}$. If for some
$\{h_n\}$ with $h_n>0$ for all n or $h_n<0$ for all n and
$\lim_{n\rightarrow\infty}h_n=0,$ $\lim_{n\rightarrow\infty}D_V(x,h_n)$
exists, then
\begin{eqnarray}
\max\big\{1-\lim_{n\rightarrow\infty}D_V(x,h_n),&& \lim_{n\rightarrow\infty}D_V(x,h_n)(p+rxI\{x\ge 0\}+\alpha xI\{x<0\})\nonumber\\
&&-(
\lambda+\delta)V(x) +\int^{x+\frac{p}{\alpha}}_0V(x-y)\dif F(y)\big \}\le0.\label{l7}
\end{eqnarray}
\end{Theorem}
\proof
We employ  a standard technique in the controlled Markov process theory, which has also been used in \cite{BenthKarlsenReikvam2001}, \cite{AlbrecherThonhauser2008} and \cite{AzcueMuler2005}.

First, we show that $V(x)$ is a viscosity super-solution  of  the HJB equation \eqref{1} on $(-\frac p\alpha,+\infty)$.
For any fixed  $x\in (-\frac p\alpha,\infty)$. Let $\phi$ be a continuously differentiable
function with $\phi(x)=V(x)$ and $V-\phi$ attaining a minimum at $x$.
For any $h>0$, define
\begin{eqnarray*}
a(x,l,h)=\left\{\begin{array}{ll}
xe^{rh}+(p-l)\int_0^h
e^{r(h-s)}\dif s,& x>0 \mbox{ and } l\geq 0 \mbox{ or } x=0  \mbox{ and } 0\leq l\leq p, \\
xe^{\alpha h}+(p-l)\int_0^h e^{\alpha(h-s)}\dif s,& -\frac p\alpha<x<0 \mbox{ and } l\geq 0 \mbox{ or } x=0 \mbox{ and } l>p.
\end{array}\right.
\end{eqnarray*}
  For any $l\geq 0$, choose an $h$ small enough such that $a(x,l,h)\in (-\frac p\alpha,0)\cup (0,\infty)$. Consider a dividend strategy $L^\prime$ that the insurer
pays out dividends continuously at rate $l$ until time $ S_1\wedge h$.
Then under the strategy $L^\prime$, ruin will not occur
before the earlier of the arrival time of the first claim $S_1$ and time $h$.
Notice that given the initial reserve $R_0=x$,  under the strategy $L^\prime$ we have
$R_t=a(x,l,t)$ for $t< S_1\wedge h$  and  $R_{S_1}=(a(x,l,S_1)-U_1)\vee (-\frac p\alpha)$ on $\{S_1\leq h\}.$ By the Dynamic Programme Principle, distinguishing two cases $S_1\geq h$ and $S_1<h$ and then conditioning on $S_1$ we have
\begin{eqnarray}
V(x)&=&\sup_{L\in \Pi} \E_{x}\left[\int_0^{S_1\wedge h}e^{-\delta s}\dif L_s+e^{-\delta (S_1\wedge h)}V(R_{S_1\wedge h})\right]\nonumber\\
  &\geq& e^{-\lambda h}\left(\int_0^h e^{-\delta s}l\dif s+e^{-\delta h}V\left(a(x,l,h)\right)\right)\nonumber\\
  &&+\int_0^h \lambda e^{-\lambda t}\dif t\Big\{\int_0^t le^{-\delta s}\dif s
   +e^{-\delta t}\int_0^{a(x,l,t)+\frac p\alpha}V\left(a(x,l,t) - u\right)\dif F(u)\Big\}.\label{l9}
\end{eqnarray}
By subtracting $V(x)$ from the last inequality  and noting that $V(x) = \phi(x) $ and $V(a(x,l,h))\geq \phi(a(x,l,h))$, we obtain
\begin{eqnarray}
  0&\geq& \frac l\delta e^{-\lambda h}(1-e^{-\delta h}) + \int_0^h
  \lambda e^{-\lambda t}\dif t\int_0^t le^{-\delta s}\dif s\nonumber\\
&&+e^{-(\lambda+\delta) h}\left[\phi\left(a(x,l,h)\right) -\phi(x)\right] + (e^{ - (\lambda + \delta)h} - 1)V(x)\nonumber\\&&
  +\int_0^h \lambda e^{-(\lambda+\delta) t}\dif t\int_0^{a(x,l,t)+\frac p\alpha}V\left(a(x,l,t) - u\right)\dif F(u).\label{l10}
\end{eqnarray}
Dividing by $h>0$ and then letting $h\downarrow 0$ yields
\begin{eqnarray}
0&\geq& l(1-\phi^\prime(x))-(\lambda+\delta)V(x)+\lambda\int_0^{x+\frac p\alpha}V(x-u)\dif F(u)\nonumber\\
&&+\left(rxI\{x>0\}\cup\{x=0,l\le p\}+\alpha xI\{x<0\}\cup\{x=0,l> p\}+p\right)\phi^\prime(x).\nonumber
\end{eqnarray}
Letting $l=0$ shows $(p+rxI\{x\ge 0\}+\alpha xI\{x<0\})\phi^\prime(x)-(\lambda+\delta)V(x)+\lambda\int_0^{x+\frac p\alpha}V(x-u)\dif F(u)\leq 0$, and  letting $l$ be large enough indicates $\phi^\prime (x)\geq 1.$

 Next, we will show that  $V(x)$ is a viscosity sub-solution  of  \eqref{1} on $(-\frac p\alpha,+\infty)$.
To this end, we employ the proof by contradiction. Assume that $V$ is not a viscosity sub-solution of \eqref{1} at some point $x$. Then we can find a constant $\eta>0$ and a continuously differentiable function $\psi_0$ with $\psi_0(x)=V(x)$, $\psi_0(y)\geq V(y)$ for all $y$ and
\begin{eqnarray}
\max\{1-\psi_0^\prime(x),\frac 1\lambda \mathcal{L}_{\psi_0}(x)\}<-2\eta.\label{12}
\end{eqnarray}
Define
\begin{eqnarray}
  \label{eq:1star}
 \psi_1(y)=\psi_0(y)+\eta \left(\frac{x-y}{x+2\frac p\alpha}\right)^2.
\end{eqnarray}
Then $\psi_1(y)$ is continuously differentiable with $  \psi_1(x)=\psi_0(x)=V(x)$  and $\psi_1^\prime(x)=\psi_0^\prime(x)$, and  by
\eqref{3} and \eqref{eq:1star} we have
$
\frac1\lambda\mathcal{L}_{\psi_1}(x)=
\frac1\lambda\mathcal{L}_{\psi_0}(x)+\int_0^{x+\frac p\alpha }\eta\left(\frac{u}{x+2\frac p\alpha}\right)^2\dif
F(u)<-2\eta+\eta=-\eta,$ which along with the fact that  $\psi_1$ is nonnegative and continuously differentiable, and
$\mathcal{L}_{\psi_1}(x)$ is continuous, indicates that there exists an $h>0$ such that
\begin{eqnarray}
\max\{1-\psi_1^\prime(y),\frac 1\lambda\mathcal{L}_{\psi_1}(y)\}<-\frac{\eta}2 \mbox{ for $y\in [x-2h,x+2h]$}.\label{15}
\end{eqnarray}

 Let $k$ be an continuously differentiable and nonnegative function with support
included in $(-1,1)$ such that $\int_{-1}^{1} k(s)\dif s=1$.
Define a function $v_n(y):\mathbb{R}\rightarrow [0,\infty)$ as the convolution
$$v_n(y)=\int_{-1}^1 \left(V(y-\frac sn)+\frac{\eta h^2}{2(x+2\frac p\alpha)^2}\right)k(s)\dif s,$$
and another function  $v:\mathbb{R}\rightarrow [0,\infty)$ by $v(y)=V(y)+\frac{\eta h^2}{2(x+2\frac p\alpha)^2}.$
Since $V(y-\frac sn)+\frac{\eta h^2}{2(x+2\frac p\alpha)^2}$  is
continuous with respect to $(y,s)$ on $[-\frac p\alpha,x+h;-1,1]$,
then we conclude that $v_n(y)$ is continuous on $[-\frac p\alpha,x+h]$.
  Moreover, $v_n(y)$ is a monotone sequence and by the dominated convergence theorem, it converges to $v(y)$. Therefore, it follows by Dini's theorem that
$v_n(y)\rightarrow v(y)$  uniformly on $[-\frac p\alpha,x+h]$.
Hence, there exists  an $n_0$ such that for all  $y\in [-\frac p\alpha,x+h]$,
\begin{eqnarray}V(y)+\frac{\eta h^2}{(x+2\frac p\alpha)^2} \geq v_{n_0}(y)\geq V(y)+\frac{\eta h^2}{4(x+2\frac p\alpha)^2}. \label{16}\end{eqnarray}

Define $f_{n_0}(y,s)=\left(V(y-\frac s{n_0})+\frac{\eta h^2}{2(x+2\frac p\alpha)^2}\right)k(s).$ It can be seen that $f_{n_0}(y,s)$ is continuous on $[-\frac p\alpha,x+h;-1,1]$.
Let $D=\{y: V(y) \mbox{
is differentiable}\}$ and $n_0(y-D)=\{n_0(y-s):s\in D\}$.
 As $V$ is
 differentiable almost everywhere,  the complement of $D$ is a null set.
Noting that
$\frac{\partial}{\partial y}f_{n_0}(y,s)=\frac{\partial}{\partial
  y}V(y-\frac s{n_0})k(s)$  on $[-\frac
p\alpha,x+h]\times n_0(y-D)$ and that
 $V^\prime(y)$, if exists, is greater than or equal
 to $1$,  it follows that for $y\in [-\frac p\alpha,x+h]$,
\begin{eqnarray}
v_{n_0}^\prime(y)=\int_{-1\leq s\leq
  1,\ s\in n_0(y-D)} \frac{\partial}{\partial y} f_{n_0}(y,s)\dif s\geq \int_{-1}^1 k(s)\dif s=1 .\label{17}
\end{eqnarray}

Construct a continuously differentiable function $\omega: \mathbb{R}\rightarrow [0,1]$ such that
$\omega(y)=1$  for $y\in [x-h,x+h]$,  $\omega(y)=0$  for $y\in(-\infty, x-2h)\cup (x+2h,\infty), $ and
 $\omega^{\prime}(y)\geq 0$  for $y\in [x-2h,x-h]$.
 Consider a function $\psi$ defined by
 \begin{eqnarray}
 \psi (y)=\omega(y)\psi_1(y)+\left(1-\omega(y)\right)v_{n_0}(y).\label{cc9}
 \end{eqnarray}
 Obviously, $\psi (x)=\psi_1(x)=V(x)$. Noting that $\psi_0\geq V$, it follows by \eqref{eq:1star}, \eqref{15}, \eqref{16}, \eqref{17} and \eqref{cc9} that
 \begin{eqnarray}
 \psi^\prime (y)
&\geq&
 \omega(y)(1+\frac \eta 2)+\omega^\prime(y)(\psi_1(y)-v_{n_0}(y))
+(1-\omega(y))\nonumber \\
&\geq&1+\omega^\prime(y)\frac{\eta}{(x+2\frac p\alpha)^2}((x-y)^2- h^2)\geq 1 \mbox{ for $y\in [-\frac p\alpha, x-h]$}. \label{9}\end{eqnarray}
 The last inequality follows by noticing  $\omega^\prime(y)\geq 0$ and $(x-y)^2- h^2\geq 0$ when $-\frac p\alpha\leq y\leq x-h$.

 By \eqref{eq:1star} and \eqref{16}, using the fact that that $V\leq \psi_0$ we obtain that for   $-\frac p\alpha\leq y-u\leq x+h$ and $h\in(0,\frac  p{2\alpha})$,
$
v_{n_0}(y-u)-\psi_1(y-u)
\leq\eta \frac{h^2}{\left(x+2\frac p\alpha\right)^2 }-\eta\left(\frac{x-y+u}{x+2\frac p\alpha }\right)^2<\frac \eta 4,$ which along with the definitions for  $\mathcal{L}$ in \eqref{3}  and  $\psi$ in \eqref{cc9}, and the fact that $\psi=\psi_1$ on $[x-h,x+h]$ indicates that  for $y\in[x-h,x+h]$,
\begin{eqnarray}
\mathcal{L}_\psi(y)
&=&\mathcal{L}_{\psi_1}(y)+\lambda\int_{y-x+h}^{y+\frac p\alpha }(1-\omega(y-u))\left(v_{n_0}(y-u)-\psi_1(y-u)\right)\dif
F(u)\nonumber\\
&<&\mathcal{L}_{\psi_1}(y)+\frac {\lambda\eta} 4.\label{cc5}
\end{eqnarray}
Write  $A=1+\frac\delta\lambda$. For any positive
 \begin{eqnarray}
 \epsilon\leq \min\{\frac{\eta h^2}{12(x+2\frac p\alpha)^2},\frac \eta{8(A-1)}\},\label{cc10}
 \end{eqnarray}
it follows by \eqref{15} and  \eqref{cc5} that
\begin{eqnarray}
\frac1\lambda\mathcal{L}_\psi(y)\leq-\frac \eta 4\leq -2(A-1)\epsilon \  \mbox{ for $y\in[x-h,x+h]$}.\label{b01}
\end{eqnarray}
From the definitions \eqref{cc9} and \eqref{eq:1star}  for $\psi$ and $\psi_1$,  respectively, by noting $0\leq w(y)\leq 1$ and using \eqref{16} it follows that for any $y$ satisfying $|y-x|\geq h$, we have
\begin{eqnarray}
\psi (y)&\geq&\omega(y)\left(\psi_0(y)+\eta \left(\frac{x-y}{x+2\frac p\alpha}\right)^2\right)+(1-\omega(y))\left(V(y)+\frac{\eta h^2}{4(x+2\frac p\alpha)^2}\right)\nonumber\\
&\geq & V(y)+3\epsilon,\label{11}
\end{eqnarray}
where the last inequality follows by the fact  $\psi_0\geq V$ and \eqref{cc10}.

From the definition \eqref{cc9} for the function $\psi$,  and the fact that all the functions $\psi$, $\omega$ and $v_{n_0}$ are continuously differentiable, we can see that $\mathcal{L}_\psi$ is  continuous. Therefore,
there exist a constant $K>0$ such that
\begin{eqnarray}
\frac 1\lambda \mathcal{L}_\psi(y)\leq K \ \mbox{ for $y\in [-\frac p\alpha ,x+h]$}.\label{012}
 \end{eqnarray}

For any fixed $\sigma$ with
\begin{eqnarray}
0<\sigma <\min\left\{\frac{\epsilon}{2\lambda K},\frac 1{4\lambda(A-1)},\frac 1\alpha\log\left(\frac{x+\frac p\alpha-\frac h2}{x+\frac p\alpha-h}\right)\right\},\label{20}\label{cc7}
\end{eqnarray}
define $\overline{\tau}=\inf\{t>0:R_t\geq x+h\}$, $\underline{\tau}=\inf\{t>0:R_t\leq x-h\}$, and  $\tau^*=\overline{\tau}\wedge (\underline{\tau}+\sigma).$
By \eqref{11} we have
\begin{eqnarray}
V(R_{\tau^*})\leq \psi(R_{\tau^*})-3\epsilon \mbox{\ \ on $\{\omega:|R_{\tau^\ast}-x|\geq h$\}.}\label{e01}
\end{eqnarray}
    Note that $R_{\bar{\tau}}=x+h$, as the surplus process has only downward jumps. Then  given the initial reserve $R_0=x$,   we have on the set $\{\omega: |R_{\tau^\ast}-x|< h\}$, \begin{eqnarray}x+h>R_{\tau^\ast}=R_{\underline{\tau}+\sigma}>x-h\geq R_{\underline{\tau}}=R_{\underline{\tau}\wedge \overline{\tau}}.\label{7a5}\end{eqnarray}
  Then it follows by \eqref{cc7}, \eqref{7a5} and noticing  that
$
  R_{t+\sigma}\leq R_te^{\alpha \sigma}+\frac p\alpha(e^{\alpha \sigma}-1)$ from the dynamics \eqref{dyn},   that given $R_0=x$,
    $$R_{\tau^\ast}=R_{\underline{\tau}+\sigma}\leq R_{\underline{\tau}}e^{\alpha \sigma}+\frac p\alpha(e^{\alpha \sigma}-1)\leq (x-h+\frac p\alpha)e^{\alpha \sigma}-\frac p\alpha< x-\frac h2\ \mbox{ on $\{\omega: |R_{\tau^\ast}-x|< h\}$},$$
    which implies that given $R_0=x$,
        $
    R_{\tau^\ast}-x<-\frac h2   \ \mbox{ on } \{\omega: |R_{\tau^\ast}-x|< h\}.$ Hence, $(x-R_{\tau^\ast})^2>\frac{h^2}4$ on the set $\{\omega:|R_{\tau^\ast}-x|< h\}$. Using this and  noticing that $\psi_0\geq V$ and that  from \eqref{cc10} we have $    \frac\eta{(x+2\frac{p}{\alpha})^2}\geq \frac{12\epsilon}{h^2}$, by the definitions \eqref{eq:1star} and \eqref{cc9}, we can show  that given the initial value $R_0=x$, $\psi(R_{\tau^\ast})=\psi_0(R_{\tau^\ast})+\eta \left(\frac{x-R_{\tau^\ast}}{x+2\frac p\alpha}\right)^2\geq V(R_{\tau^\ast})+3\epsilon$  on $\{\omega:|R_{\tau^\ast}-x|< h\}$.
This along with \eqref{e01} shows
\begin{eqnarray}
V(R_{\tau^*})\leq \psi(R_{\tau^*})-3\epsilon \ \mbox{ given $R_0=x$}.\label{21}
\end{eqnarray}

Note that from \eqref{9} we have
 $ {\psi}^\prime(R_t)\geq 1$ for $R_t\in [-\frac p\alpha,x-h]$ and  that from \eqref{15} and \eqref{cc9} we have $ {\psi}^\prime(R_t)= {\psi}_1^\prime(R_t)\geq 1$ for $[x-h,x+h]$.
 Then by noticing that $R_t\in [-\frac p\alpha,x+h]$ for $t\in [0,\tau^\ast]$,  we conclude that ${\psi}^\prime(R_t)\geq 1$  for $t\in [0,\tau^\ast]$.
 Then by setting $\phi$ and $\tau$ in \eqref{26} to be $\psi$ and $\tau^*$, respectively,  we have
 \begin{eqnarray}
 \psi(R_{\tau^*})e^{-\delta \tau^*}- \psi(R_0)\leq \int_0^{\tau^*}\mathcal{L}_\psi(R_{t-})e^{-\delta t}\dif t-\int_0^{\tau^*} e^{-\delta t}\dif L_t+\mathcal{M}_{\tau^*}(\psi).\label{7a10}
\end{eqnarray}

Note that given $R_0=x$, we have $R_{t-}\in [x-h,x+h]$ for $t\in[0,\bar{\tau}\wedge\underline{\tau}]$ and $R_{t-}\in [-\frac p\alpha,x+h]$ for $t\in [\bar{\tau}\wedge\underline{\tau},\tau^\ast]$.
From \eqref{b01} and  \eqref{012}, it follows that
\begin{eqnarray}
\int_0^{\tau^*}\mathcal{L}_\psi(R_{t-})e^{-\delta t}\dif t&=&\int_0^{\overline{\tau}\wedge \underline{\tau}}\mathcal{L}_\psi(R_{t-})e^{-\delta t}\dif t+\int_{\overline{\tau}\wedge \underline{\tau}}^{\tau^*}\mathcal{L}_\psi(R_{t-})e^{-\delta t}\dif t\nonumber\\
&\leq& -(A-1)2\epsilon \lambda\int_0^{\overline{\tau}\wedge \underline{\tau}}e^{-\delta t}\dif t+\lambda K\int_{\overline{\tau}\wedge \underline{\tau}}^{\tau^*}e^{-\delta t}\dif t\nonumber\\
 &=&-(A-1)2\epsilon \lambda\int_0^{\tau^*}e^{-\delta t}\dif t+\lambda((A-1)2\epsilon+ K)\int_{\overline{\tau}\wedge \underline{\tau}}^{\tau^*}e^{-\delta t}\dif t\nonumber\\
 &\leq&-(A-1)2\epsilon \lambda\int_0^{\tau^*}e^{-\delta t}\dif t+\lambda((A-1)2\epsilon+ K)\sigma\nonumber\\
 &<&-(A-1)2\epsilon \lambda\int_0^{\tau^*}e^{-\delta t}\dif t+\epsilon,\label{cc11}
\end{eqnarray}
where the second last inequality follows by noticing $\tau^\ast-\bar{\tau}\wedge\underline{\tau}\leq\sigma$ and
the last inequality follows by \eqref{20}.

Given the initial reserve $R_0=x$, it follows from \eqref{21}, \eqref{7a10} and \eqref{cc11}  that
\begin{eqnarray}
V(R_{\tau^*})e^{-\delta \tau^*}
&<&\psi(R_{\tau^*})e^{-\delta \tau^*}-2\epsilon e^{-\delta \tau^*}\nonumber\\
&=&(\psi(R_{\tau^*})e^{-\delta \tau^*}-\psi(x))+(\psi(x)-2\epsilon e^{-\delta \tau^*})\nonumber\\
&\leq&-(A-1)2\epsilon\lambda
\int_0^{\tau^*} e^{-\delta t}\dif
t-\int_0^{\tau^*}e^{-\delta s}\dif L_t+\mathcal{M}_{\tau^*}(\psi)\nonumber\\
&&+(\psi(x)-2\epsilon
e^{-\delta \tau^*})+\epsilon.\label{029}\end{eqnarray}
 As $ \int_0^{\tau^*}e^{-\delta s}\dif s=\frac{1-e^{-\delta \tau^*}}\delta$ and $A=1+\frac \delta\lambda$, from \eqref{029} we obtain
\begin{eqnarray}
&&V(R_{\tau^*})e^{-\delta \tau^*}+\int_0^{\tau^*}e^{-\delta s}\dif L_t
\leq \mathcal{M}_{\tau^*}(\psi)+\psi(x)-\epsilon .\label{new7}
\end{eqnarray}
Noting that $\mathcal{M}_{t}(\psi)$ is a super-martingale with zero-expectation, we have $\E[\mathcal{M}_{\tau^*}(\psi)]\leq 0$.
As a result, taking conditional expectation on \eqref{new7}  yields $
V(x)=\sup_{L\in\Pi}\E_x[
\int_0^{\tau^*}e^{-\delta s}\dif L_t+V(R_{\tau^*})e^{-\delta \tau^*}]\leq \psi(x)-\epsilon,$
which contradicts the fact $V(x)=\psi(x)$.

\noindent (ii) By Theorem \ref{thm2.2}, it follows that for any $h_n$ with $\lim_{n\rightarrow\infty}h_n=0$, \begin{eqnarray}
\lim_{n\rightarrow\infty}D_V(x,h_n)\geq 1.\label{l12}
\end{eqnarray}
  Consider a sequence $h_n^\prime$ with $h_n^\prime\downarrow 0$ as $n\rightarrow \infty$ such that $\lim_{n\rightarrow \infty}\frac{V(a(x,l,h_n^\prime))-V(x)}{h_n^\prime}$ exists.
 Following the same lines as in the proof for (i), it can be shown that  $\eqref{l10}$ also holds when $h$ and $\phi(\cdot)$ there being replaced by $h_n^\prime$ and $V(\cdot)$, respectively. Dividing both sides of  the newly obtained inequality by $h_n^\prime$ and then letting $n\rightarrow \infty$ yields for $l\ge 0$,
 \begin{eqnarray*}
&&0\geq l(1-\lim_{n\rightarrow\infty}D_V(x,a(x,l,h_n^\prime)-x))-(\lambda+\delta)V(x)+\lambda\int_0^{x+\frac p\alpha}V(x-u)\dif F(u)+\big[p+\nonumber\\
&&(rI\{x>0\}\cup\{x=0, l\le p\}+\alpha I\{x<0\}\cup\{x=0, l> p\})x\big]\lim_{n\rightarrow\infty}D_V(x,a(x,l,h_n^\prime)-x).
\end{eqnarray*}
By letting $l=0$ we have \begin{eqnarray}&&(p+rxI\{x\ge 0\}+\alpha x(\{x<0\})\lim_{n\rightarrow\infty}D_V(x,a(x,0,h_n^\prime)-x)-(\lambda+\delta)V(x)\nonumber\\
&&+\lambda\int_0^{x+\frac p\alpha}V(x-u)\dif F(u)\leq 0.\label{l11}
\end{eqnarray}
For any
$\{h_n\}$ with $h_n>0$ such that
$\lim_{n\rightarrow\infty}h_n=0,$ and $\lim_{n\rightarrow\infty}D_V(x,h_n)$
exists, we can find a subsequence $\{h_{n_k}\}\subset \{a(x,l,h_n^\prime)-x\}$.
Therefore, $\lim_{n\rightarrow\infty}D_V(x,h_n)
=\lim_{k\rightarrow\infty}D_V(x,h_{n_k})
=\lim_{n\rightarrow\infty}D_V(x,a(x,0,h_n^\prime)-x).$
It follows by \eqref{l12} and \eqref{l11} that
\begin{eqnarray}
\max\big\{1-\lim_{n\rightarrow\infty}D_V(x,h_n),&& \lim_{n\rightarrow\infty}D_V(x,h_n)(p+rxI\{x\ge 0\}+\alpha xI\{x<0\})\nonumber\\
&&-(
\lambda+\delta)V(x) +\int^{x+\frac{p}{\alpha}}_0V(x-y)\dif F(y)\big \}\le0.\label{l17}
\end{eqnarray}

For any sequence
$\{h_n\}$ with $h_n<0$ such that
$\lim_{n\rightarrow\infty}h_n=0,$ and $\lim_{n\rightarrow\infty}D_V(x,h_n)$
exists, by repeating the above argument by replacing all $x$ there by $x-c(x,l,h)$ (i.e., conditioning on the initial reserve $R_0=x-c(x,l,h)$), where
$$c(x,l,h)=\left\{\begin{array}{ll}
x(1-e^{-rh})+(p-l)\int_0^he^{-rs)}\dif s& x>0 \mbox{ and } l\geq 0 \mbox{ or } x=0  \mbox{ and } 0\leq l\leq p \\
x(1-e^{-\alpha h})+(p-l)\int_0^he^{-\alpha s)}\dif s&-\frac p\alpha<x<0 \mbox{ and } l\geq 0 \mbox{ or } x=0 \mbox{ and } l>p,
\end{array}\right.$$
and noticing that $a(x-c(x,l,h),l,h)=x$, we can show that \eqref{l17} is also true.  \hfill $\square$\\

Next we will show that the value function $V(x)$  is the unique nonnegative, nondecreasing and locally Lipschitz continuous viscosity solution of \eqref{1} satisfying a linear growth condition and the boundary condition $V(-\frac p\alpha)=0$. We start with the following comparison principle.
\begin{Lemma}\label{lemma32}
Let $\overline{u}(x)$ and $\underline{u}(x)$ be a  nonnegative viscosity super-solution and sub-solution, respectively. Assume that for both $u=\overline{u}(x)$ and $\underline{u}(x)$, the function $u$ is continuous on $[-\frac p\alpha, \infty)$ and locally Lipschitz continuous on $(-\frac p\alpha,\infty)$, and satisfies  $u(-\frac p\alpha)=0$ and    $u(x)\leq c_1x+c_2$ for some constants $c_1$ and $c_2$ . Then  $\underline{u}(x)\leq \overline{u}(x)$ for all $x\geq -\frac p\alpha$.\label{T1}
\end{Lemma}

From Theorem \ref{thm2.1}, Theorem \ref{thm2.3} and Theorem \ref{thm3.1}, we know that the value function $V(x)$ is a nondecreasing and nonnegative viscosity solution of the HJB equation \eqref{1} that is  locally Lipschitz continuous on $(-\frac p\alpha,\infty)$,  satisfies a linear growth condition, and  fulfills the boundary condition $V(-\frac p\alpha)=0$. Consider any other viscosity solution $W(x)$ of \eqref{1} that fulfils the same conditions. Since $V(x)$ is also a super-solution and   and $W(x)$ is also a sub-solution, by Lemma \ref{T1} we conclude that $V(x)\geq W(x)$ for all $x\geq -\frac p\alpha$. This leads to the following theorem stating the uniqueness of the value function as a viscosity solution of \eqref{1}.
\begin{Theorem}
The value function $V(x)$ is the unique nondecreasing and nonnegative viscosity solution of the HJB equation \eqref{1} that

i) is locally Lipschitz on $(-\frac p\alpha,\infty)$,

ii) satisfies a linear growth condition, and

iii) fulfills the boundary condition $V(-\frac p\alpha)=0$.
\end{Theorem}

As an immediate result of Lemma \ref{T1}, we arrive at the Verification Theorem as follows.

\begin{Theorem}
For any strategy $L\in \Pi$, if $V_L$ is an locally Lipshcitz continuous
viscosity super-solution of HJB equation \eqref{1}, then $V_L=V$, i.e. $L$ is an optimal dividend strategy.
\end{Theorem}
{\proof} Obviously, $V_L$ is nonnegative and nondecreasing and $V_L(-\frac p\alpha)=0$.
Since $ V_L\leq V$, it is true that $V_L$ also satisfies the linear growth condition. Therefore, by Lemma \ref{T1} we know that $V_L\geq V$. Consequently, $V_L=V$.
\hfill $\square$\\

\begin{Lemma}\label{lemma32}
  Let $\overline{u}(x)$ and $\underline{u}(x)$ be a viscosity
  super-solution and sub-solution of the HJB equation \eqref{1} on
  $[b_0,\infty)$, respectively. Assume that for both
  $u=\overline{u}(x)$ and $\underline{u}(x)$, the function $u$ is
  continuous on $[b_0, \infty)$ and satisfies $u(x)\leq c_1x+c_2$ for
  some constants $c_1$ and $c_2$. If $ \underline{u}(b_0)\leq
  \overline{u}(b_0)$, then $\underline{u}(x)\leq \overline{u}(x)$ for
  all $x\geq b_0$.\label{T1}
\end{Lemma}
The proof is in Appendix.

 \begin{Remark}\label{ap30} By Lemma \ref{lemma32} it is obvious that
   for any given constant $c$, there is at most one viscosity
   solution, $u$, of the equation \eqref{1} on $[b_0,\infty)$ that
   satisfies the initial condition $u(b_0)=c$ and the linear growth
   condition.
\end{Remark}

\begin{Lemma}\label{T2} Let $\Pi_{\overline{x}}$ be the set of admissible strategies such that the controlled reserve $R_t$ is less than or equal to $\overline{x}$ for all $t>0$.
If for some $\bar{x}>0$, $\overline{u}(x)$ is a    nonnegative, nondecreasing and locally  Lipshcitz continuous super-solution of the HJB equation \eqref{1} on $(-\frac p\alpha, \bar{x})$, then $\overline{u}(x)\geq  \sup_{L\in\Pi_{\overline{x}}}V_L(x)$ for all $x\in[-\frac p\alpha, \overline{x})$.
\end{Lemma}
\proof
i) We can prove this by showing that for any dividend strategy $L\in \Pi_{\overline{x}}$, $V_L(x)\leq \overline{u}(x)$ for $x\in[-\frac p\alpha, \overline{x})$. For any continuous super-solution $\overline{u}$ of the HJB equation \eqref{1} on $(-\frac p\alpha, \overline{x})$, consider a function $\overline{v}(x)$ with   $\overline{v}(x)=  0 $ for $x<-\frac p\alpha$, $\overline{v}(x)=\overline{u}(x)$ for $x\in [-\frac p\alpha,\overline{x})$ and $\overline{v}(x)=                   \overline{u}(\overline{x})$ for $x\geq \overline{x}$.
                    Consider a sequence of nonnegative functions $v_n(x)=\int_{-\infty}^\infty \overline{v}(x-y)n\phi(ny)\dif y$ for $x\in [-\frac p\alpha, \overline{x}],$
where $\phi(x)$ is a nonnegative, even and continuously differentiable function with its support included in $(-1,1)$ such that $\int_{-1}^1 \phi(x)\dif s=1$.
It can be seen that $v_n(x)$ is nonnegative and nondecreasing, and satisfies
\begin{eqnarray}
v_n(x)\leq \overline{u}(\overline{x}), \mbox{ for } x\in [-\frac p\alpha, \overline{x}].\label{b5}
\end{eqnarray}
Using the standard techniques in  real analysis (eg \cite{WheedenZygmund1977}), we can show that $v_n$ is continuously differentiable on $[-\frac p\alpha, \overline{x}]$,
\begin{eqnarray}
 &&\mbox{$v_n(x)$ converges to $\overline{u}(x)$ uniformly on  $[-\frac p\alpha, \overline{x}]$; and} \label{xx1}\\
   &&\mbox{$v_n^\prime(x)$ converges to $\overline{u}^\prime(x)$ almost everywhere.\label{xx2}}
  \end{eqnarray}
Noting from Definition \ref{Def} (ii), $
1\leq \overline{u}^\prime (x)\leq \frac{\lambda+\delta}{p+rxI\{x\geq 0\}+\alpha xI\{x<0\}}\overline{u}(x)
$ for $x$ such that $\overline{u}^\prime(x)$ exists, we can obtain
\begin{eqnarray}
&&1\leq v_n^\prime (x)\ \mbox{ for } x\in[-\frac p\alpha, \overline{x}].\label{b6}
\end{eqnarray}

From now on in the proof of this lemma, we assume $ x\in[-\frac p\alpha, \overline{x})$.
  By setting $u$ and $\tau$ in \eqref{26} to be $\overline{v}_n$ and $t\wedge \tau$, respectively,  and  then taking expectation, we obtain
\begin{eqnarray}
\E_x[v_n(R_{t\wedge  T})e^{-\delta (t\wedge  T)}]
\leq v_n(x)+\E_x\left[\int_0^{t\wedge T} e^{-\delta s}\mathcal{L}_{v_n}(R_{s-})\dif s\right]-\E_x\left[\int_0^{t\wedge T}e^{-\delta s}\dif L_s\right].\label{b7}
\end{eqnarray}
Notice that under any strategy $L\in\Pi_{\overline{x}}$ the controlled reserve is below or at $\overline{x}$. Then by \eqref{b5} we have $\int_0^{R_{s-}+\frac p\alpha} v_n(R_{s-}-y)\dif
F(y)\leq \overline{u}(\overline{x})$. Further note that $p+rR_{s-}I\{R_{s-}\geq 0\}+\alpha R_{s-}I\{R_{s-}<0\}\ge 0$ for $x\in[0,T]$.  Hence,  by using \eqref{3}, \eqref{b5}, \eqref{b6},  the monotone convergence and the dominated convergence we can obtain
\begin{eqnarray}
\lim_{t\rightarrow\infty}\E_x[\int_0^{ T\wedge t} e^{-\delta s}\mathcal{L}_{v_n}(R_{s-})\dif s]&=&\E_x\left[\int_0^{ T} e^{-\delta s}\mathcal{L}_{v_n}(R_{s-})\dif s\right].\label{0b6}
\end{eqnarray}
Letting $t\rightarrow\infty$ on both sides of  \eqref{b7} and then using \eqref{0b6}, the dominated convergence and the monotone convergence  yields
\begin{eqnarray}
0\leq v_n(x)+\E_x\left[\int_0^{ T} e^{-\delta s}\mathcal{L}_{v_n}(R_{s-})\dif s\right]-V_L(x).\label{b8}
\end{eqnarray}
Since under any strategy $L \in\Pi_{\overline{x}}$,
 the controlled reserve $R_{s-}\leq \overline{x}$,
  by \eqref{xx1} and \eqref{xx2} it can be shown that $
  \mathcal{L}_{v_n}(R_{s-})\rightarrow \mathcal{L}_{\overline{u}}(R_{s-}) \ \ a.e.-\pr.$
  Using this and Fatou's Lemma, taking  $\limsup_{n\rightarrow \infty}$ on \eqref{b8} yields
$
0\le \overline{u}(x)+\E_x\left[\int_0^{ T} e^{-\delta s}\mathcal{L}_{\overline{u}}(R_{s-})\dif s\right]-V_L(x).
$ As $\overline{u}$ is a super-solution, so  $\mathcal{L}_{\overline{u}}(R_{s-})\leq 0$ $a.e.-\pr$. Then it follows that $V_L(x)\leq \overline{u}(x).$
Consequently, $V(x)=\sup_{L\in\Pi_{\overline{x}}}V_L(x)\leq \overline{u}(x).$
\hfill $\square$\\

\bigskip
Define an operator $\mathcal{G}$ by
\begin{eqnarray}
\mathcal{G}_u(x)&=&p+rxI\{x\geq 0\}+\alpha x I\{x< 0\}-(\lambda+\delta)u(x)\nonumber\\
&&+\lambda\int_0^{x+\frac p\alpha }u(x-y)\dif
F(y). \label{7a21}
\end{eqnarray}
\begin{Theorem}\label{T3}
If for some $\bar{x}\in(-\frac p\alpha,\infty)$, $\mathcal{G}_V(\bar{x})=0$, then
$V(x)=\sup_{L\in \Pi_{\bar{x}}}V_L(x)$ on $[-\frac p\alpha,\bar{x}]$,
where $\Pi_{\overline{x}}$ is defined same as in Lemma \ref{T2}.
\end{Theorem}
\begin{Theorem}\label{T4} Let $\Pi_{\bar{x}}$ be defined same as in Lemma \ref{T2}. If there exists an $\bar{x}\in (-\frac p\alpha,\infty)$ such that $V^\prime(\bar{x})=1$, then $
V(x)=\sup_{L\in \Pi_{\bar{x}}}V_{L}({x}) \ \mbox{ for } x\in[-\frac p\alpha, \overline{x}].$
\end{Theorem}

As an immediate consequence of Theorem \ref{T2}, Theorem \ref{T3} and Theorem \ref{T4}, we obtain the following theorem.
\begin{Theorem}\label{T5}
If for some $\bar{x}>-\frac p\alpha$ with either $\mathcal{G}_{V}(\bar{x})=0$ or $V^\prime(\bar{x})=1$,  then for any nonnegative, nondecreasing and locally Lipschitz continuous  super-solution $\overline{u}(x)$  of the HJB equation \eqref{1} on $(-\frac p\alpha,\bar{x}]$ which satisfies $\overline{u}(x)\leq c_1+c_2x$ for some constants $c_1$ and $c_2$, and the boundary condition $\overline{u}(-\frac p\alpha)=0$,  we have $\overline{u}(x)\geq V(x)$ on $(-\frac p\alpha,\bar{x}]$. Furthermore, if for some strategy $L\in\Pi_{\bar{x}}$ the function $V_L$ is an absolutely continuous super-solution to the HJB equation \eqref{1} on $(-\frac p\alpha,\bar{x}]$, then  $V(x)=V_{L}({x})$ for all $x\in(-\frac p\alpha,\bar{x}]$.
\end{Theorem}

For any $y\geq -\frac p\alpha$, define $G_y(x)=V(x)$ if $ x\leq y$ and $G_y(x)=V(y)+x-y$ if $x>y$.
\begin{Theorem}\label{theorem39} i) If $G_y$ is a super-solution to the HJB equation \eqref{1} on $(y,\infty)$, then $G_y=V$ on $[-\frac p\alpha,\infty)$.\\
ii) If for some $\bar{x}>-\frac p\alpha$ with either $\mathcal{G}_{V}(\bar{x})=0$ or $V^\prime(\bar{x})=1$, and for some $y<\bar{x}$, $G_y$ is a super-solution of the HJB equation \eqref{1} on $(y,\bar{x}]$, then  $G_y(x)= V(x)$ on $[-\frac p\alpha,\bar{x}]$.
\end{Theorem}
\proof  First we show that $G_y$ is a viscosity super-solution to the HJB equation \eqref{1} on $(-\frac p\alpha,y]$. For any fixed $x\in[-\frac p\alpha,y]$, let $\phi$ be any continuously differentiable function with $\phi(x)=G_y(x)$ and $G_y-\phi$ reaches minimum at $x$.
Then by Remark \ref{cc50} ii) we obtain
\begin{eqnarray}
\limsup_{h\uparrow 0}\frac{G_y(x)-G_y(x-h)}{h}\leq \lim_{n\rightarrow\infty}D_V(x,a(x,h_n^\prime)-x)\leq \liminf_{h\downarrow 0}\frac{G_y(x+h)-G_y(x)}{h}.\label{cc61}
\end{eqnarray}
Notice that
$
\limsup_{h\uparrow 0}\frac{G_y(x)-G_y(x-h)}{h}=\limsup_{h\uparrow 0}\frac{V(x)-V(x-h)}{h}$
and
that $\liminf_{h\downarrow 0}\frac{V_y(x+h)-V_y(x)}{h}$ equals $\liminf_{h\downarrow 0}\frac{G_y(x+h)-G_y(x)}{h}$ if $x\in[-\frac p\alpha,y)$ and equals $1$ if $x=y$.
As a result, using \eqref{cc61} and $\liminf_{h\downarrow 0}\frac{V_y(x+h)-V_y(x)}{h}\geq 1$ yields
$
\limsup_{h\downarrow 0}\frac{V(x)-V(x-h)}{h}\leq \phi^\prime(x)\leq \liminf_{h\downarrow 0}\frac{V(x+h)-V(x)}{h},
$
which by Remark \ref{cc50} ii) again implies that $V-\phi$ reaches minimum at $x$. Since $V$ is a viscosity super-solution of \eqref{1}, we have $\max\{1-\phi^\prime(x),\mathcal{L}_{V,\phi}(x)\}\leq 0.$
Hence, by noticing $\mathcal{L}_{V,\phi}(x)=\mathcal{L}_{G_y,\phi}(x)$ for $x\in(-\frac p\alpha,y]$, we have $\max\{1-\phi^\prime(x),\mathcal{L}_{G_y,\phi}(x)\}\leq 0\mbox{ for $x\in(-\frac p\alpha,y]$}. $
Consequently, $G_y$ is a viscosity super-solution on $(-\frac p\alpha,y].$

\noindent i) If $G_y$ is a viscosity super-solution on $(y,\infty),$ then it is a super-solution on $(-\frac p\alpha,\infty)$. Also note that $G_y$ satisfies the linear growth condition. Then by Theorem \ref{T2} i), we have $G_y\geq V$ on $(-\frac p\alpha,\infty)$. Noticing  that $G_y\leq V$, therefore, $G_y=V$ on $[-\frac p\alpha,\infty)$.

\noindent ii) If $G_y$ is a viscosity super-solution on $(y,\bar{x}],$ then it is a super-solution on $(-\frac p\alpha,\bar{x}]$. By Theorem \ref{T5}, we have $G_y\geq V$ on $(-\frac p\alpha,\bar{x}]$. Noticing by definition that $G_y\leq V$, therefore, $G_y=V$ on $(-\frac p\alpha,\bar{x}]$.
\hfill $\square$

\section{The optimal dividend strategy}
In this section, we show that there exists an optimal dividend strategy and the optimal  strategy is a band strategy, that is, the optimal strategy at any time is to pay no dividends, pay out at a rate same as the premium incoming rate or a positive lump sum, depending on the current reserve at that time. We also show that under certain condition, when the reserve is negative, the optimal strategy is to pay no dividends.

We start with the following definition for three sets.
\begin{Definition} Define
$\mathcal{A}=\{x\in [-\frac p\alpha,\infty):\mathcal{G}_V(x)=0\},$
$\mathcal{B}=\{x\in [-\frac p\alpha,\infty):V^\prime(x)=1 \mbox{ and } \mathcal{G}_V(x)<0\},$ and
$\mathcal{C}=(\mathcal{A} \cup \mathcal{\beta})^c.$
\end{Definition}
The sets defined above will play a crucial role in proving the existence of and characterizing the optimal dividend strategy.
  we can prove the following useful properties of these sets.
\begin{Lemma} \label{L0c}The following properties hold.\\
(a) $\mathcal{A}$ is nonempty and closed.\\
 (b) $\mathcal{B}$ is nonempty and left-open. And there exists a $y$ such that $(y,\infty)\subset \mathcal{B}.$\\
 (c) If $(x_0,x_1]\subset\mathcal{B}$ and $x_0 \notin \mathcal{B}$, then $x_0\in \mathcal{A}$.\\
 (d) $\mathcal{C}$ is right-open.
\end{Lemma}

Based on the above three sets and their characteristics, we define the following dividend strategy, which will be shown to be the optimal one.
\begin{Definition} \label{def44}Let $L^\ast$ be a dividend strategy defined as follows:

(a) If $R^{L^\ast}_{t-}\in\mathcal{A} $, the insurer pays out dividends at the same rate as the premium incoming rate, i.e.

 $\dif L^\ast_t=\left(p+rR_{t-}^{L^\ast}I\{R_{t-}^{L^\ast}\geq 0\}+\alpha R_{t-}^{L^\ast}I\{R_{t-}^{L^\ast}< 0\}\right)\dif t$ if $R^{L^\ast}_{t-}\in \mathcal{A} $.

 (b) If $R^{L^\ast}_{t-}\in\mathcal{B} $, then by Lemma \ref{L0c} (c) there exists an $x_0\in \mathcal{A} $ with $x_0<R^{L^\ast}_{t-}$ such that $(x_0, R^{L^\ast}_{t-}]\subset \mathcal{B} $. At time $t$, the insurer pays out a lump sum $R^{L^\ast}_{t-}-x_0$ as dividends, ie

  $L^\ast_t- L^\ast_{t-}=R^{L^\ast}_{t-}-x_0$ if $R^{L^\ast}_{t-}\in \mathcal{B} $, where   $x_0=\inf\{x:(x, R^{L^\ast}_{t-}]\subset \mathcal{B} \}$.

  (c) If $R^{L^\ast}_{t-}\in\mathcal{C} $, then the insurer pays out no dividends at the moment.
\end{Definition}

In the following, we prove that the strategy $L^\ast$ constructed above is an optimal dividend strategy.
\begin{Theorem}\label{thm42}
The strategy $L^*$ defined in Definition \ref{def44} is optimal, i.e. $V(x )=V_{L^*}(x )$ for all $x\geq -\frac p\alpha$.
\end{Theorem}
\proof By Lemma \ref{L0c} it follows that there exists some $\underline{x} =\inf\{x:(x,\infty)\subset\mathcal{B} \}$.

 Let $\mathcal{H}$ be a set of continuous functions $f: [-\frac p\alpha,\infty)\rightarrow [0,\infty)$ with $f(x )=x-\underline{x} +f(\underline{x}  )$ for $x>\underline{x} $.

Define the distance $\rho(f_1,f_2)=\max_{x\geq -\frac p\alpha}|f_1(x)-f_2(x)|$ for  $f_1,f_2\in\mathcal{H}$.

Define an operator $\mathcal{T}$ as follows:
\begin{eqnarray}
\mathcal{T}_f(x )=\E_{x}\left[\int_{0}^{S_1}e^{-\delta s}\dif L^*_s + e^{-\delta S_1}f(R^{L^\ast}_{S_1} )\right].\label{Toperator}
\end{eqnarray}
 Noting that for any $x\geq \underline{x} $, we have $(\underline{x} ,\infty)\subset\mathcal{B} $ and $\underline{x} \in\mathcal{A} $ , by  using Definition \ref{def44}(b) with  $x_0=\underline{x}$ and \eqref{Toperator} we get
\begin{eqnarray}
\mathcal{T}_f(x )
=x-\underline{x} +\mathcal{T}_f(\underline{x}  )\ \mbox{ for $x\ge \underline{x}$.}\label{setb}\end{eqnarray}
As a result, $\mathcal{T}_f\in \mathcal{H}$ for any $f\in \mathcal{H}.$ 
Note that
\begin{eqnarray*}
|T_{f_1}(x )-T_{f_2}(x )|&=&|\E_{x} [e^{-\delta S_1}(f_1(R^{L^\ast}_{S_1} )-f_2(R^{L^\ast}_{S_1} ))]|\nonumber\\
&\leq&\frac{\lambda }{\lambda+\delta}\rho(f_1,f_2),\label{contraction}
\end{eqnarray*}
where the last inequality follows by the fact that $S_1$ is an exponential random variable with mean $\frac1{\lambda}$.
Therefore,  $\mathcal{T}$ is a contraction on $\mathcal{H}$ and thus has a unique fixed point in $\mathcal{H}$.

According to the structure of $L^*$ (Definition \ref{def44}), we can see that the process $L^*$ is a Markov process and therefore the controlled reserve process under $L^*$ is also a Markov process. By the Markov property  and \eqref{Toperator}, it is obvious that $V_{L^*}$ is a fixed point of $\mathcal{T}$ in $\mathcal{H}$.
So to prove $V=V_{L^*}$ it is sufficient to show that $V\in\mathcal{H}$ and $V$ is also a fixed point of $\mathcal{T}$.

Obviously, $V\in[0,\infty)$.  Moreover, since $(\underline{x} ,\infty)\subset\mathcal{B} $, then $V^\prime(x )=1$ for all $x> \underline{x} $. As a result, $V(x )=V(\underline{x}  )+x-\underline{x} $ for all $x\geq \underline{x} $. Consequently, we can conclude that $V\in\mathcal{H}$.

 Assume $x\in \mathcal{A} $. By  the definition of $L^*$ , we can see that given $R_0=x$, $\dif L^\ast_t=(p+rxI\{x\geq0\}+\alpha xI\{x<0\})\dif t$ for all time $t$ before the arrival $S_1$ of the first claim. Therefore, by \eqref{Toperator} we obtain that
\begin{eqnarray}
\mathcal{T}_V(x )&=&\E_{x}\big[\int_{0}^{S_1}(p+rxI\{x\geq 0\}+\alpha xI\{x<0\})e^{-\delta s}\dif s+e^{-\delta S_1}V(x-U_{1} )]\nonumber\\
&=&\frac{(p+rxI\{x\geq 0\}+\alpha x I\{x<0\})}{\lambda+\delta}+\int_{0}^{\infty}\lambda e^{-\lambda  t}e^{-\delta t} \dif t \int_{0}^{x+\frac p\alpha}V(x-y )\dif F (y)\nonumber\\
&=& \frac{(p+rxI\{x\geq 0\}+\alpha xI\{x<0\})+\lambda \int_{0}^{x+\frac p\alpha}V(x-y )\dif F (y)}{\lambda+\delta} \mbox{ for $x\in \mathcal{A} $}.\label{aset}\end{eqnarray}
    It follows by \eqref{aset} and the equality $\mathcal{G}_V(x )=0$ for $x\in\mathcal{A}$ that
\begin{eqnarray}
\mathcal{T}_V(x )
=V(x )\mbox{ for $x\in \mathcal{A} $}. \label{aequal}
\end{eqnarray}

 For any $x\in\mathcal{B} $, we can find an $x_0<x$ such that $(x_0,x]\subset \mathcal{B} $ and $x_0\in \mathcal{A} $, which implies $V^\prime (y )=1$ for $y\in(x_0,x]$. Therefore, $V(x )=x-x_0+V(x_0 )$. By the definition of $L^*$, we know that a lump sum of $x-x_0$ will be paid out as dividends immediately. Then it follows from \eqref{Toperator} and \eqref{aequal} that
\begin{eqnarray}
\mathcal{T}_V(x )=x-x_0+\mathcal{T}_V(x_0 )=x-x_0+V(x_0 )=V(x ) \mbox{ for $x\in \mathcal{B} $}.\label{bequal}
\end{eqnarray}

Now we look at the case $x\in\mathcal{C} $. Since $\mathcal{C}$ is right open, there exists an $x_1$ such that $(x,x_1)\subset \mathcal{C}$ and $x_1\notin \mathcal{C}$. As $\mathcal{B}$ is left open, so $x_1 \in\mathcal{A}$. By the definition for $L^*$ we know that given the initial reserve $R_0=x$, the insurance company pays out no dividends until the reserve reaches $x_1$ or the arrival ($S_1$) of the first claim.
Consider a function $a(\cdot)$ which satisfies $\dif a(t)=\left(p+ra(t)I\{a(t)\geq 0\}+\alpha a(t)I\{a(t)< 0\}\right)\dif t$ and $a(0)=x$. Recall that $t_0(x,x_1)$ is the time it will take for this dynamics to reach $x_1$. It can be seen that given $R_0=x$, $R_t=a(t)$ for all $x<S_1\wedge t_0(x,x_1)$,  and  $R_{S_1}=a(S_1)-U_{1}$ if $S_1<t_0(x,x_1)$.

 By Markov property it follows that for any $t\geq 0$,
\begin{eqnarray}
\mathcal{T}_f(x )&=&\E_{x}\left[\int_{0}^{S_1}e^{-\delta s}\dif L^*_s+e^{-\delta S_1}f(R^{L^\ast}_{S_1} );S_1\leq t\right]\nonumber\\
&&+\E_{x}\left[\int_{0}^{t}e^{-\delta s}\dif L^*_s + e^{-\delta t}\mathcal{T}_f(R^{L^\ast}_{t} );S_1> t\right].\label{new4}
\end{eqnarray}

By setting $t$ and $f$ in \eqref{new4} by $t_0(x,x_1)$ and $V$, respectively, and by noting $\mathcal{T}_V(x_1)=V(x_1)$ because $x_1\in\mathcal{A}$, it follows that
\begin{eqnarray}
\mathcal{T}_V(x )&=&\E_{x}\left[e^{-\delta S_1}V(a(S_1)-U_{1} )I\{S_1\leq t_0(x,x_1)\}+e^{-\delta t_0(x,x_1)}\mathcal{T}_V(x_1)I\{S_1>t_0(x,x_1)\}\right]\nonumber\\
&=&\int_{0}^{t_0(x,x_1)}\lambda e^{-\lambda  t}e^{-\delta t} \dif t \int_{0}^{a(t)+\frac p\alpha}V(a(t)-y )\dif F (y)+  e^{-(\lambda+\delta) t_0(x,x_1)}V(x_1).\label{JEaa20}
\end{eqnarray}
Let $\mathcal{D} =\{x>0: V'(x ) \mbox{ exists }\}$  and $t\in \overline{\mathcal{D}} := \{y: a(y)\in \mathcal{D} \}$.
As $V(x)$ is differentiable almost everywhere, the Lebesgue measure of $\mathcal{D}^c $ is $0$.
Noting that $V(a(t))$ is differentiable for $a(t)\in \mathcal{D} $,  the complement of $\overline{\mathcal{D}} $ has a zero Lebesgue measure, too.

\noindent Notice that for any $y$ such that $V^\prime(y )$ exists we have \begin{eqnarray}
&&(p+ryI\{y\geq 0\}+\alpha y I\{y<0\})V^\prime(y )-(\lambda+\delta)V(y )+\lambda \int_{0}^{y+\frac p\alpha}V(y-z )\dif F (z)=0.\nonumber\\\label{0aa92}
\end{eqnarray}
It follows from \eqref{JEaa20} and \eqref{0aa92} that
\begin{eqnarray}
\mathcal{T}_V(x )&=&\int_{\overline{\mathcal{D}} \cap (0,t_0(x,x_1))} e^{-(\lambda+\delta)t}\Bigg((\lambda+\delta)V(a(t) )\nonumber\\
&&
-\Big(p+ra(t)I\{a(t)\geq 0\}+\alpha a(t)I\{a(t)<0\}\Big) V^\prime(a(t) )\Bigg)\dif t\nonumber\\
&&+  e^{-(\lambda+\delta) t_0(x,x_1)}V(x_1)\nonumber\\
&=&\int_{\overline{\mathcal{D}} \cap (0,t_0(x,x_1))} \dif \big( e^{-(\lambda+\delta)t}V(a(t) \big)+  e^{-(\lambda+\delta) t_0(x,x_1)}V(x_1)\nonumber\\
&=& V(x )-e^{-(\lambda+\delta) t_0(x,x_1)}V(x_1)+e^{-(\lambda+\delta) t_0(x,x_1)}V(x_1)\nonumber\\
&=&V(x),\ \mbox{ for $x\in\mathcal{C} $.}\label{cequal}
\end{eqnarray}
Combining \eqref{aequal}, \eqref{bequal} and \eqref{cequal} shows that $V(\cdot)$ is a fixed point of $\mathcal{T}$. This completes the proof.
\hfill $\square$\\

Now we have shown that like the Cram\'er-Lundberg cases respectively with and without interest, the optimal strategy is also a band strategy in the absolute ruin case. Intuitively, we would think that under the optimal strategy, there should be no dividends if the company is in deficit. In the following we  will prove this rigorously.
\begin{Lemma}\label{thml1} For any fixed $x_0\in(-\frac{p}{\alpha},\infty)$, there
exists a unique in $(x_0,\infty)$ differentiable, strictly
increasing and positive solution $u$ on $[x_0,\infty)$ to
the equation
\begin{eqnarray}
0&=&(p+rxI\{x\ge0\}+\alpha xI\{x< 0\})u^\prime(x)-(\lambda+\delta)u(x)
\nonumber \\&&+\lambda\int^{x-x_0}_0u(x-y)\dif F(y)
+\lambda\int^{x+\frac{p}{\alpha}}_{x-x_0}V(x-y)\dif F(y) \label{l6}
\end{eqnarray} with boundary condition $u(x_0)=V(x_0).$
\end{Lemma}
\proof  First we show that there is a such solution on $[x_0,x_0+h)$
for $h=\frac{p+\alpha x_0I\{x_0<0\}}{2(2\lambda+\delta)}$.\\
 Let
$\mathcal{H}[x_0,x_0+h)$ denote the set of continuous, increasing and
positive functions on $[x_0,x_0+h)$. Define an operator $\mathcal{T}$ that for any
$u\in \mathcal{H}[x_0,x_0+h)$,
\begin{eqnarray}\mathcal{T}_u(x)=\int^{x}_{x_0}
\frac{(\lambda+\delta)u(s)-\lambda\int^{s-x_0}_0u(s-y)\dif F(y)
-\lambda\int^{s+\frac{p}{\alpha}}_{s-x_0}V(s-y)\dif F(y)}
{p+rsI\{s\ge 0\}+\alpha sI\{s<0\}}\dif s +V(x_0).\label{l5}
\end{eqnarray}
We will show that $\mathcal{T}$ is a contraction
on $\mathcal{H}[x_0,x_0+h)$.\\
For any $u\in \mathcal{H}[x_0,x_0+h)$, as both $u$ and $V$ are increasing and $u(x_0)=V(x_0)$, we get
\begin{eqnarray*}
&&\lambda\int^{x-x_0}_{0}u(x-y)\dif F(y)+\lambda\int^{x+\frac{p}{\alpha}}_{x-x_0}V(x-y)\dif F(y)\\
&\leq&
\lambda\int^{x-x_0}_{0}u(x)\dif F(y) +
\lambda\int^{x+\frac{p}{\alpha}}_{x-x_0}u(x)\dif F(y)\le\lambda u(x).
\end{eqnarray*}
Define $||u||=\sup_{x\in[x_0,x_0+h)}|u(x)|$.\\
It follows by \eqref{l5} that for any $x\in[x_0, x_0+h)$ and
$u_1,u_2\in\mathcal{H}[x_0,x_0+h)$,
\begin{eqnarray*}|\mathcal{T}_{u_1}(x)-\mathcal{T}_{u_2}(x)|&\le& \int^{x}_{x_0}\frac{(\lambda+\delta) ||u_1-u_2|| ^{x_0+h}_{x_0}+\lambda ||u_1-u_2|| ^{x_0+h}_{x_0}}
{p+rsI\{s\ge0\}+\alpha sI\{s<0\}}\dif s\\
&\le&
\frac{(2\lambda+\delta)h||u_1-u_2||^{x_0+h}_{x_0}}{p+\alpha
x_0I\{x_0<0\}}\leq\frac12 ||u_1-u_2||^{x_0+h}_{x_0}.
\end{eqnarray*}
Therefore, $\mathcal{T}$ is a contraction on
$\mathcal{H}[x_0,x_0+\varepsilon)$. As a result, there exists a
unique $u\in \mathcal{H}[x_0,x_0+\varepsilon)$ such that
$u(x)=\mathcal{T}_u(x)$, i.e.,
$$u(x)=\int^{x}_{x_0}\frac{(\lambda+\delta)u(s)-\lambda\int^{s-x_0}_{0}u(s-y)\dif F(y)
-\lambda\int^{s+\frac{p}{\alpha}}_{s-x_0}V(s-y)\dif F(y)}
{p+rsI\{s\ge0\}+\alpha sI\{s<0\}}\dif s+V(x_0),$$
which implies
\begin{eqnarray*}
u^\prime(x)=\frac{(\lambda+\delta)u(x)-\lambda\int^{x-x_0}_{0}u(x-y)\dif F(y)
-\lambda\int^{x+\frac{p}{\alpha}}_{x-x_0}V(x-y)\dif F(y)}
{p+rxI\{x\ge0\}+\alpha xI\{x<0\}}\\
\mbox{ for $x\in [x_0,x_0+h)$.}
\end{eqnarray*}
This completes the  proof of the existence and uniqueness of an
positive,increasing and differentiable solution to (6)
on $[x_0,x_0+h)$.

Similarly, we can prove the existence and uniqueness of a solution
to (6) on $[x_0+h,x_0+2h)$ fulfilling the required properties. Repeating
the above process, we can prove the existence of a unique solution
to (6) on $[x_0,\infty)$, which is differentiable, increasing and
positive.
\hfill $\square$\\

\begin{Theorem}\label{thml2}
 (i) For any $x\in\mathcal{A}$, $V(x)$ is
differentiable and $V^\prime (x)=1$.\\
(ii) For any $(x_0,x_1)\subset \mathcal{C}$, $V(x)$ is
differentiable on $(x_0,x_1)$, and $V^\prime (x)>1$ for
$x\in(x_0,x_1)$.
\end{Theorem}
\proof
(i) Consider any $x\in \mathcal{A}$. Choose
sequences $h_n^+>0$ and $h_n^-<0$ with
$\lim_{n\rightarrow\infty}h_n^\pm=0$, such that
$\lim_{n\rightarrow\infty}\frac{V(x+h_n^+)-V(x)}{h_n^+}=\limsup_{h\downarrow
0}\frac{V(x+h)-V(x)}{h}$ and
$\lim_{n\rightarrow\infty}\frac{V(x+h_n^-)-V(x)}{h_n^-}=\limsup_{h\uparrow
0}\frac{V(x+h)-V(x)}{h}.$
As
$x\in\mathcal{A}$, $\mathcal{G}_V(x)=0$. Then
it follows by \eqref{3} and Theorem \ref{thm3.1} (ii)  that
$\lim_{n\rightarrow\infty}D_V(x,h_n^{\pm})\le 1,$
which implies
$$\limsup_{h\uparrow 0}\frac{V(x+h)-V(x)}{h}\le 1 \ \mbox{ and }\ \limsup_{h\downarrow 0}\frac{V(x+h)-V(x)}{h}\le 1.$$
As $\liminf_{h\rightarrow 0}\frac{V(x+h)-V(x)}{h}\ge 1,$ we conclude that
$\lim_{h\rightarrow 0}\frac{V(x+h)-V(x)}{h}= 1.$

\noindent (ii) Use proof by contradiction.
%
Note that for any $x\in\mathcal{C}\cap (-\frac p\alpha,0)$, if $V(x)$ is
differentiable, then
$0=(p+rxI\{x\geq 0\}+\alpha xI\{x< 0\})V^\prime(x)-(\lambda+\delta)V(x)+\lambda\int^{x+\frac{p}{\alpha}}_{0}V(x-y)\dif F(y),$
which can be rewritten as
\begin{eqnarray}0=&&(p+rxI\{x\geq 0\}+\alpha xI\{x< 0\})V^\prime(x)-(\lambda+\delta)V(x)+\lambda\int^{x-x_0}_{0}V(x-y)\dif F(y)\nonumber\\
&&+
\lambda\int^{x+\frac{p}{\alpha}}_{x-x_0}V(x-y)\dif F(y).\label{l221}
\end{eqnarray}
Then by \eqref{l221} and Lemma \ref{thml1}, we conclude that $V(x)$ is equal to the unique
solution of \eqref{l6} on $(x_0,x_0+h)$ and therefore is differentiable on
$(x_0,x_1)$.

By Theorem \ref{thm2.2}, we know that for any $x\in\mathcal{C}$, if $V^\prime(x)$ exists, then $V^\prime(x)\geq 1$. By the definition of the set $\mathcal{C}$, we know that $V^\prime(x)$, if exists, can not be $1$. If $V^\prime(x)=1$, then $x$ belongs to either $\mathcal{A}$ or $\mathcal{B}$.  Therefore,
  $V^\prime(x)\neq 1$ for all $x\in(x_0,x_1)$.
 \hfill $\square$\\

\begin{Theorem}\label{thml3}
 Assume $\alpha > \lambda+\delta$. The following statements hold.

(i) $\mathcal{A}\cap (-\frac{p}{\alpha},0)$ consists of
isolated points only.

  (ii) For any $x_0\in
\mathcal{A}\cap(-\frac{p}{\alpha},0)$, we can find an $h>0$ such that
$(x_0,x_0+h)\subset \mathcal{B}$.

(iii) $\mathcal{B}\cap(-\frac{p}{\alpha},0)=\emptyset$.

\end{Theorem}

\proof
 Consider any $x_1$ and $x_2$ with $-\frac{p}{\alpha}<x_1<x_2<0$  and  $[x_1,x_2)\subset
\mathcal{A}$ such that  $V(x)$ is differentiable on $[x_1,x_2)$,  $V^\prime(x_1)=1$ and
\begin{eqnarray}(p+\alpha
x)V^\prime(x)=(\lambda+\delta)V(x)-\lambda\int^{x+\frac{p}{\alpha}}_0V(x-y)\dif F(y) \mbox{ for $x\in[x_1,x_2)$.}\label{1222}
\end{eqnarray}
By setting $x$ in the above equality to be $x_1$ and $x_1+\epsilon$, respectively, and using the newly obtained equations, we can obtain that for any $\varepsilon\in(0,x_2-x_1)$,
\begin{eqnarray} &&\frac{(p+\alpha
x_1)(V^\prime(x_1+\varepsilon)-V^\prime(x_1))}{\varepsilon}\nonumber\\
&=&-\alpha
V^\prime(x_1+\varepsilon)+(\lambda+\delta)\frac{V(x_1+\varepsilon)-V(x_1)}{\varepsilon}-\lambda
I(x_1,\varepsilon),\label{fm001}\end{eqnarray}
where $
I(x_1,\varepsilon)=\frac{\int^{x_1+\varepsilon+\frac{p}{\alpha}}_0V(x_1+\varepsilon-y)\dif F(y)
-\int^{x_1+\frac{p}{\alpha}}_0V(x_1-y)\dif F(y)}{\varepsilon}.$\\
By noticing that $V^\prime(x_1)=1, V^\prime(x_1+\varepsilon)\ge 1$, $I(x_1,\varepsilon)\ge 0$ and $\lambda+\delta<\alpha$, from \eqref{fm001} we obtain
$(p+\alpha x_1)\frac{V^\prime(x_1+\varepsilon)-V^\prime(x_1)}{\varepsilon}
\le -\alpha+(\lambda+\delta)(1+\frac{o(\varepsilon)}{\varepsilon})<0$
 for  small $\varepsilon$.
 As a result,
\begin{eqnarray}V^\prime(x_1+\varepsilon)<V^\prime(x_1)=1 \ \mbox{ for   $\varepsilon$ ($\varepsilon>0$) small enough}.\label{lstar}\end{eqnarray}

\noindent (i) Use proof by contradiction. Assume that $x_0\in
\mathcal{A}\cap (-\frac{p}{\alpha},0)$ and it is not isolated. Then, as
$\mathcal{A}$ is closed, we can find an $h>0$ such that either
$[x_0,x_0+h]\subset \mathcal{A}$ or $[x_0,x_0+h]\subset
\mathcal{A}$. Use $[x_1,x_2]$ to denote $[x_0-h,x_0]$ if $[x_0-h,x_0]\subset
\mathcal{A}$, and $[x_0-h,x_0]$, otherwise. Then $[x_1,x_2]\subset \mathcal{A}$.\\
It follows by Theorem \ref{thml2} (i) that
\begin{eqnarray}V^\prime(x)=1\  \mbox{ for $x\in [x_1,x_2]$}. \label{l20}
  \end{eqnarray}
  Therefore, according to the definition for $\mathcal{A}$, we have $\mathcal{L}_V(x)=\mathcal{G}_V(x)=0$  for all $x\in [x_1,x_2]$, which is equivalent to $(p+\alpha
x)V^\prime(x)=(\lambda+\delta)V(x)-\lambda\int^{x+\frac{p}{\alpha}}_0V(x-y)\dif F(y)$
for $x\in[x_1,x_2]$.  Then  by \eqref{lstar} it follows that
$V^\prime(x_1+\varepsilon)<V^\prime(x_1)=1$ for  small positive $\varepsilon$,
 which is a contradiction
to \eqref{l20}.

\noindent (ii)
Assume that there exists an $x_0\in\mathcal{A}\cap(-\frac p\alpha,0)$, such that we can find an $h>0$ satisfying $(x_0,x_0+h)\nsubseteq\mathcal{B}$. Then
$(x_0,x_0+h)\subset \mathcal{C}$, because $\mathcal{A}$ consists of isolated points only and both $\mathcal{B}$ and
$\mathcal{C}$ are half open. Hence, it follows by Theorem \ref{thml2} (ii) that $V(x)$ is differentiable on
$(x_0,x_0+h)$ and $V^\prime(x)>1$ for
$x\in(x_0,x_0+h)$. Hence, $V$ is a solution to the HJB equation $\eqref{1}$ and therefore, \eqref{1222} holds for $x\in(x_0,x_0+h)$.
As $x_0\in\mathcal{A}$, we have $V^\prime(x_0)=1$, which along with the definition for $\mathcal{A}$ implies that \eqref{1222} also holds for
$x=x_0$.
Then by setting $x_1$ and $x_2$ in \eqref{lstar} as $x_0$ and $x_0+h$, respectively, it follows that
$V^\prime(x_0+\varepsilon)<V^\prime(x_0)=1$  for  small positive $\varepsilon$,
 which is a contradiction
to the fact that $V^\prime(x_0+\varepsilon)\ge 1$.

\noindent (iii)
Assume $\mathcal{B}\cap(-\frac{p}{\alpha},0)\neq\phi$. Then there exist
$x_0$ and $x_1$, such that $-\frac p\alpha<x_0<x_1<0$, $[x_0,x_1)\subset \mathcal{B}$
and $x_0\in\mathcal{A}$. Therefore, by Theorem \ref{thml2} (i) and the definition for $\mathcal{B},$ we get $V^\prime(x)=1$ for
$x\in[x_0,x_1)$, which implies
$V(x)=x-x_0+V(x_0) \mbox{ for $x\in[x_0,x_1)$}.$
Note that $\mathcal{G}_V(x_0)=0$. Then for $x\in(x_0,x_1)$,
\begin{eqnarray}\mathcal{G}_V(x)&=&\mathcal{G}_V(x)-\mathcal{G}_V(x_0)\nonumber\\
&=&\alpha(x-x_0)-(\lambda+\delta)(x-x_0)
+\lambda\int^{x+\frac{p}{\alpha}}_0V(x-y)\dif F(y)
-\lambda\int^{x_0+\frac{p}{\alpha}}_0V(x_0-y)\dif F(y)\nonumber\\
&>&0, \label{l2star}
\end{eqnarray}
where the last inequality follows by $\alpha>\lambda+\delta$ and the
fact that $V$ is nonnegative and increasing.
Since
$x\in(x_0,x_1)\subset\mathcal{B}$, we have $\mathcal{G}_V(x)<0$, which
contradicts the inequality \eqref{l2star}.
\hfil $\square$\\

\begin{Theorem}\label{thml4} If $\alpha>\lambda+\delta$, $(-\frac{p}{\alpha},0)\subset \mathcal{C}$.
\end{Theorem}
\proof By Theorem \ref{thml3} (iii), it follows that \begin{eqnarray}(-\frac{p}{\alpha},0)\cap
 \mathcal{B}=\emptyset.\label{l22}
\end{eqnarray}
So it is sufficient to show that $(-\frac{p}{\alpha},0)\cap
\mathcal{A}=\emptyset$. If this is not true, then we can find an $x_0\in
\mathcal{A}\cap(-\frac{p}{\alpha},0)$. By Theorem \ref{thml3}, it follows that there exist an $h>0$ such that
$(x_0,x_0+h)\subset \mathcal{B}$, which contradicts  \eqref{l22} by noting $x_0+h\in(-\frac{p}{\alpha},0)$ for small $h>0$.
\hfill $\square$

\begin{Remark}  Theorem \ref{thm42} and Theorem \ref{thml4} together imply that  if $\alpha>\lambda+\delta$, under  the optimal strategy $L^\ast$ the company   will pay no dividends when the reserve is negative. In other words, if $\alpha>\lambda+\delta$ it is optimal to pay no dividends when the reserve is negative.
\end{Remark}

\section{Conclusion}
 We studied the dividend optimization problem of an insurance corporation, of which the surplus is modeled by a compound Poisson model with credit and debit interest. The company earns interest when the reserve is positive, and can refinance to settle its claims when the reserve is negative but above the critical level.  The company controls the dividend pay-out dynamically and seeks to maximize the expected total discounted dividends until ruin. We proved that the value function is the unique viscosity solution satisfying certain conditions of the associated Hamilton-Jacob-Bellman equation, that the optimal strategy is a band strategy, and that it is optimal to pay no dividends when the reserve is negative. This result provides theoretical justification to the regulation of no dividend payments when the surplus is in deficit.\\

\noindent{\bf\large Acknowledgements}\\
 I would like to thank Feng Chen for many valuable comments and suggestions and the referee for advices on improving the presentation of the paper. Financial support by Australian School of Business Research Grants, University of New South Wales, is gratefully acknowledged.

\numberwithin{equation}{section}
 \renewcommand{\theequation}{A-\arabic{equation}}
  \setcounter{equation}{0}  
  \section*{APPENDIX}  

\noindent In this appendix, we present the proofs to  Lemma \ref{lemma32}, Theorem \ref{T3}, Theorem \ref{T4} and Lemma \ref{L0c} \\

\bigskip

\noindent {\bf Proof of Lemma \ref{lemma32}}
We employ a proof by contradiction.

Assume that there exists a $x_0\in(-\frac p\alpha,\infty)$ such that
$\underline{u}(x_0)> \overline{u}(x_0)$.

For any constant $\gamma>0$, define functions for $x\geq -\frac p\alpha$,
$$\overline{u}_{\gamma}(x)=e^{-\gamma x}\overline{u}(x)\ \ \mbox{ and }\ \ \underline{u}_{\gamma}(x)=e^{-\gamma x}\underline{u}(x).$$
By the fact that both the functions  $\overline{u}$ and $\underline{u}$  are locally Lipschitz continuous and  bounded by a linear function, it can be easily shown that $\overline{u}_{\gamma}(x)$ and $\underline{u}_{\gamma}(x)$ are both bounded and Lipschitz continuous on $(-\frac p\alpha,\infty)$, too, which implies that  there exists some constant $m>0$ such that
\begin{eqnarray}
\frac{\overline{u}_{\gamma}(y)-\overline{u}_{\gamma}(x)}{y-x}\leq m \ \ \mbox{ and }\
\frac{\underline{u}_{\gamma}(y)-\underline{u}_{\gamma}(x)}{y-x}\leq m   \label{010}\ \ \mbox{ for } x,y\in(-\frac p\alpha,\infty).
\end{eqnarray}

For $\rho>0$, consider a function $\phi_\rho: [-\frac p\alpha,\infty)\times [-\frac p\alpha,\infty)\rightarrow \mathbb{R}$ given by \begin{eqnarray}\phi_\rho(x,y)=\underline{u}_\gamma(x)-\overline{u}_{\gamma}(y)-\frac \rho 2 (x-y)^2-\frac{2m}{\rho^2(y-x)^2+\rho}.\label{050}
\end{eqnarray}

Note that we can find a $\gamma_1>0$ such that $\underline{u}_\gamma(x_0)-\overline{u}_\gamma(x_0)> 0$ for all $\gamma\in (0,\gamma_1]$, and that $\underline{u}(-\frac p\alpha)=\overline{u}(-\frac p\alpha)=0$ and $\lim_{x\rightarrow \infty}\underline{u}_\gamma(x)=\lim_{x\rightarrow \infty}\overline{u}_\gamma(x)=0$. Then we can define \begin{eqnarray}
 M=\max_{x\geq -\frac p\alpha}(\underline{u}_\gamma(x)-\overline{u}_\gamma(x))\ \mbox{ and }\
M_\rho=\max_{x,y\geq -\frac p\alpha}\phi_\rho(x,y).\label{cc15}
\end{eqnarray}
 Then $0<M<\infty$ and $M$ has a maximizer denoted by
 $x^*$, and $M_\rho$ also has a maximizer, denoted by $(x_\rho,y_\rho)$ here.

\noindent
 Noting that
\begin{eqnarray}
M_\rho\geq \phi_\rho(x^*,x^*)=M-\frac{2m}{\rho},\label{3080}
\end{eqnarray}
then it follows that
\begin{eqnarray}
\liminf_{\rho\rightarrow \infty}M_\rho\geq M>0. \label{060}
\end{eqnarray}

 Let $(\rho_n)_{n\in \mathbb{N}}$ be a sequence  tending to $\infty$ as $n\rightarrow \infty$ such that $(x_{\rho_n},y_{\rho_n})$ converges as $n\rightarrow \infty$. Use $(\bar{x},\bar{y})$ to denote the  limit of $(x_{\rho_n},y_{\rho_n})$ as $n\rightarrow \infty$. We will show in the following that
  \begin{eqnarray}
 \bar{x}=\bar{y}.\label{a4}
 \end{eqnarray}
 If this is not true, then $|\bar{x}-\bar{y}|>0$.  By noticing
\begin{eqnarray}M_{\rho_n}=
\underline{u}_\gamma(x_{\rho_n})-\overline{u}_{\gamma}(y_{\rho_n})-\frac {\rho_n}2 (x_{\rho_n}-y_{\rho_n})^2-\frac{2m}{\rho_n^2(y_{\rho_n}-x_{\rho_n})^2
+\rho_n},\label{3084}
\end{eqnarray}$\lim_{n\rightarrow \infty}\underline{u}_\gamma(x_{\rho_n})=\underline{u}_\gamma(\bar{x})$ and $\lim_{n\rightarrow \infty}\overline{u}_\gamma(y_{\rho_n})=\overline{u}_\gamma(\bar{y})$, we obtain
\begin{eqnarray*}\lim_{n\rightarrow\infty}M_{\rho_n}=
\underline{u}_\gamma(\bar{x})-\overline{u}_{\gamma}(\bar{y})-\frac {\lim_{n\rightarrow\infty}\rho_n}2 (\bar{x}-\bar{y})^2=-\infty,
\end{eqnarray*}
which is a contradiction to \eqref{060}.

Next, we  show that for any constant $\hat{x}\geq -\frac p\alpha$ with $\underline{u}_\gamma(\hat{x})\leq \overline{u}_{\gamma}(\hat{x})$,
\begin{eqnarray}
\bar{x}=\bar{y}\neq \hat{x}.\label{apm01}
\end{eqnarray}
We use proof by contradiction again.  Suppose  $
\bar{x}=\bar{y}=\hat{x}$. Then  for any $\epsilon^\prime>0$ we can find an $\delta^\prime>0$ such that $\underline{u}_\gamma(x)-\overline{u}_{\gamma}(x)< \epsilon$ for all $x$ satisfying $ |x-\hat{x}|<\delta^\prime$. Note that $\lim_{n\rightarrow\infty}x_{\rho_n}=\lim_{n\rightarrow\infty}y_{\rho_n}
=\bar{x}=\hat{x}$. Hence, there exists an $N^\prime>0$ such that for all $n\ge N^\prime$,  $|x_{\rho_n}-\hat{x}|<\delta^\prime$ and $|y_{\rho_n}-\hat{x}|<\delta^\prime$. Therefore, $M_{\rho_n}=\phi_{\rho_n}(x_{\rho_n},y_{\rho_n})\le \underline{u}_\gamma(x_{\rho_n})- \overline{u}_{\gamma}(y_{\rho_n})\le \epsilon^\prime$. Consequently,
$\limsup_{n\rightarrow\infty}M_{\rho_n}\le0$, which contradicts \eqref{060}.

  Noting that $\underline{u}(-\frac p\alpha)=\overline{u}(-\frac p\alpha)=0$. By \eqref{apm01}, we can conclude immediately that
 \begin{eqnarray}
\bar{x}=\bar{y}\neq \hat{x}.\label{a2}
\end{eqnarray}

By observing that
\begin{eqnarray*}\lim_{y\rightarrow \infty} \phi_\rho(x,y)=\lim_{y\rightarrow \infty}\left(\underline{u}_\gamma(x)-\overline{u}_{\gamma}(y)-\frac \rho 2 (x-y)^2-\frac{2m}{\rho^2(y-x)^2+\rho}\right)=-\infty,
\end{eqnarray*}
 we conclude that \begin{eqnarray}
 \bar{y}<\infty\label{a3}.
 \end{eqnarray}
 Combining \eqref{a4}, \eqref{a2} and \eqref{a3} yields $\bar{x}=\bar{y}\in (-\frac p\alpha,\infty).$

  As $x_{\rho_n}$ and $y_{\rho_n}$ converge to $\bar{x}$ and $\bar{y}$, respectively,  we can find an $N_1$ such that for all $n\geq N_1$,  \begin{eqnarray}&&x_{\rho_n},\ y_{\rho_n}\in (-\frac p\alpha,\infty).\label{a6}
 \end{eqnarray}

 Now we introduce two more functions
 $$\zeta_{\rho}(x)=\overline{u}_\gamma(x_\rho)+\frac{\rho}2 (x-y_\rho)^2
 +\frac{2m}{\rho^2(y_\rho-x)^2+\rho}+\phi_\rho(x_\rho,y_\rho),$$
 and
 $$\varphi_{\rho}(y)=\underline{u}_\gamma(y_\rho)-\frac{\rho}2 (x_\rho-y)^2
 -\frac{2m}{\rho^2(x_\rho-y)^2+\rho}-\phi_\rho(x_\rho,y_\rho).$$
 It can be easily shown that for all $n\geq N_1$, $\zeta_{\rho_n}$ and $\varphi_{\rho_n}$ are both continuously differentiable.
 Furthermore, $\underline{u}_\gamma(x)-\zeta_{{\rho_n}}(x)=\phi_{\rho_n}(x,y_{\rho_n})-\phi_{\rho_n}(x_{\rho_n},y_{\rho_n})$ attains its maximum $0$ at $x_{\rho_n}$, and
 $\overline{u}_\gamma(y)-\varphi_{{\rho_n}}(y)=-\phi_{\rho_n}(x_{\rho_n},y)+\phi_{\rho_n}(x_{\rho_n},y_{\rho_n})$ reaches its minimum $0$ at $y_{\rho_n}$.
 Since $\underline{u}$ and $\overline{u}$ are respectively viscosity sub and super-solutions of \eqref{1}, by the definition for viscosity solutions we can see that $\underline{u}_\gamma$ and $\overline{u}_\gamma$ are respectively viscosity sub and super-solutions of the following equation
 \begin{eqnarray*}
&&\max\big\{1-e^{\gamma x}(\gamma u(x)+u^\prime (x)),\left(p+rxI\{x\geq 0\}+\alpha x I\{x< 0\}\right)\times
\nonumber\\
&&(\gamma u(x)+u^\prime(x))-(\lambda+\delta)u(x)+\lambda\int_0^{x+\frac p\alpha }u(x-y)e^{-\gamma y}\dif
F(y)\big\}=0.
\end{eqnarray*}
Therefore, by Definition \ref{alt} we can obtain  that for $n\geq N_1$,
\begin{eqnarray}
&&\max\big\{1-e^{\gamma x_{\rho_n}}(\gamma \underline{u}_\gamma(x_{\rho_n})+\zeta_{{\rho_n}}^\prime (x_{\rho_n})),\left(p+rx_{\rho_n}I\{x_{\rho_n}\geq 0\}+\alpha x_{\rho_n} I\{x_{\rho_n}< 0\}\right)\times
\nonumber\\
&&(\gamma \underline{u}_\gamma(x_{\rho_n})+\zeta_{{\rho_n}}^\prime(x_{\rho_n}))
-(\lambda+\delta)\underline{u}_\gamma(x_{\rho_n})+\lambda\int_0^{x_{\rho_n}+\frac p\alpha }\underline{u}_\gamma(x_{\rho_n}-y)e^{-\gamma y}\dif
F(y)\big\}\geq 0,\nonumber\\\label{001}
\end{eqnarray}
and
\begin{eqnarray}
&&\max\big\{1-e^{\gamma y_{\rho_n}}(\gamma \overline{u}_\gamma(y_{\rho_n})+\varphi_{{\rho_n}}^\prime (y_{\rho_n})),\left(p+ry_{\rho_n}I\{y_{\rho_n}\geq 0\}+\alpha y_{\rho_n} I\{y_{\rho_n}< 0\}\right)\times
\nonumber\\
&&(\gamma \overline{u}_\gamma(y_{\rho_n})+\varphi_{{\rho_n}}^\prime(y_{\rho_n}))
-(\lambda+\delta)\overline{u}_\gamma(y_{\rho_n})+\lambda\int_0^{y_{\rho_n}+\frac p\alpha }\overline{u}_\gamma(y_{\rho_n}-y)e^{-\gamma y}\dif
F(y)\big\}\leq 0.\nonumber\\
\label{002}
\end{eqnarray}
Use $B_1$, $B_2$ to represent  the first and second terms in the curly brackets on the left-hand side of \eqref{001}, respectively, and $D_1$, $D_2$ to represent the first and second terms on the left-hand side of \eqref{002}, respectively. Then $\max\{B_1,B_2\}\geq 0\geq \max\{D_1,
D_2\}.$  So at least one of the inequalities $B_1\geq D_1$ and $B_2\geq D_2$ holds.

(i) First, assume that $B_2\geq D_2$ is true. Noticing that
 \begin{eqnarray}
 \zeta_{{\rho_n}}^\prime(x_{\rho_n})=
 \varphi_{{\rho_n}}^\prime(y_{\rho_n})={\rho_n} (x_{\rho_n}-y_{\rho_n})
 +\frac{4m(y_{\rho_n}-x_{\rho_n})}{\left(\rho_n(y_{\rho_n}-x_{\rho_n})^2+1\right)^2},
 \label{a12}
 \end{eqnarray}
   by substitutions for $\zeta_{{\rho_n}}^\prime(x_{\rho_n})$ and $ \varphi_{{\rho_n}}^\prime(y_{\rho_n})$ by \eqref{a12}, it follows immediately that
\begin{eqnarray}
&&\left(p+ry_{\rho_n} I\{y_{\rho_n}\geq 0\}+\alpha y_{\rho_n} I\{y_{\rho_n}< 0\}\right)\times\nonumber\\
&&\left(\gamma \overline{u}_\gamma(y_{\rho_n})+{\rho_n} (x_{\rho_n}-y_{\rho_n})
 +\frac{4m(y_{\rho_n}-x_{\rho_n})}{\left(\rho_n(y_{\rho_n}-x_{\rho_n})^2+1\right)^2}\right)\nonumber\\
 &&-\left(p+rx_{\rho_n} I\{x_{\rho_n}\geq 0\}+\alpha x_{\rho_n} I\{x_{\rho_n}< 0\}\right)\times\nonumber\\
&&\left(\gamma \underline{u}_\gamma(x_{\rho_n})+{\rho_n} (x_{\rho_n}-y_{\rho_n})
 +\frac{4m(y_{\rho_n}-x_{\rho_n})}{\left(\rho_n(y_{\rho_n}-x_{\rho_n})^2+1\right)^2}\right)\nonumber\\
&&+(\lambda+\delta)\left(\underline{u}_\gamma(x_{\rho_n})
-\overline{u}_\gamma(y_{\rho_n})\right)\nonumber\\
&\leq& \lambda\left(\int_0^{x_{\rho_n}+\frac p\alpha }\underline{u}_\gamma(x_{\rho_n}-y)e^{-\gamma y}\dif
F(y)-\int_0^{y_{\rho_n}+\frac p\alpha }\overline{u}_\gamma(y_{\rho_n}-y)e^{-\gamma y}\dif
F(y)\right).\label{011}
\end{eqnarray}
Notice that $\phi_{\rho_n}(x_{\rho_n},x_{\rho_n})+\phi_{\rho_n}(y_{\rho_n},y_{\rho_n})\leq 2\phi_{\rho_n}(x_{\rho_n},y_{\rho_n})$, i.e.
\begin{eqnarray*}&&\underline{u}_\gamma(x_{\rho_n})-
\overline{u}_{\gamma}(x_{\rho_n})+\underline{u}_\gamma(y_{\rho_n})-
\overline{u}_{\gamma}(y_{\rho_n})-\frac{4m}{{\rho_n}}\\
&&
\leq 2\left(\underline{u}_\gamma(x_{\rho_n})-\overline{u}_{\gamma}(y_{\rho_n})-\frac {\rho_n} 2 (x_{\rho_n}-y_{\rho_n})^2-\frac{2m}{{\rho_n}^2(y_{\rho_n}-x_{\rho_n})^2+{\rho_n}}\right).\label{051}
\end{eqnarray*}
Rearranging terms gives
\begin{eqnarray*}
{\rho_n} (x_{\rho_n}-y_{\rho_n})^2 &\leq& \underline{u}_\gamma(x_{\rho_n})-\underline{u}_\gamma(y_{\rho_n})
+\overline{u}_{\gamma}(x_{\rho_n})-\overline{u}_{\gamma}(y_{\rho_n})+ \frac{4m(y_{\rho_n}-x_{\rho_n})^2}{{\rho_n}(y_{\rho_n}-x_{\rho_n})^2+1} \nonumber\\
&\leq & 2m|y_{\rho_n}-x_{\rho_n}|+4m(y_{\rho_n}-x_{\rho_n})^2,
\end{eqnarray*}
where the last inequality follows by \eqref{010}.  As a result,
\begin{eqnarray*}
|y_{\rho_n}-x_{\rho_n}|\leq \frac{2m}{{\rho_n}-4m} \mbox{ for ${\rho_n}>4m$}. \label{0012}
\end{eqnarray*}

As $\underline{u}_\gamma$ and $\overline{u}_\gamma$ are both bounded,  taking limits $\lim_{n\rightarrow\infty}$ on $\eqref{011}$ yields
\begin{eqnarray}
&&\gamma\left(p+r\bar{x} I\{\bar{x}\geq 0\}+\alpha \bar{x} I\{\bar{x}< 0\}\right) \left( \overline{u}_\gamma(\bar{x})
  - \underline{u}_\gamma(\bar{x})\right)\nonumber\\&&
+(\lambda+\delta)(\underline{u}_\gamma(\bar{x})
-\overline{u}_\gamma(\bar{x}))\nonumber\\
&\leq& \lambda\left(\int_0^{\bar{x}+\frac p\alpha }\left(\underline{u}_\gamma(\bar{x}-y)-\overline{u}_\gamma(\bar{x}-y)\right)
e^{-\gamma y}\dif
F(y)\right)\ \mbox{ for $\gamma>0$}\label{ab1}\\
&\leq& \lambda M\ \mbox{for  $\gamma\in (0,\gamma_1)$},\label{cc16}
\end{eqnarray}
where the last inequality follows from \eqref{cc15}.

By choosing $\gamma<\min\left\{\frac{\delta}{2(p+r\bar{x} I\{\bar{x}\geq 0\})},\gamma_1\right\},$ it follows immediately from \eqref{cc16} that
\begin{eqnarray}
\underline{u}_\gamma(\bar{x})-\overline{u}_\gamma(\bar{x})
<\frac{\lambda}{\lambda+\frac{\delta}2}M<M.\label{061}\end{eqnarray}
 On the other hand, from \eqref{060}  we get  $
M\leq \liminf_{\rho\rightarrow \infty}M_\rho\leq \lim_{n\rightarrow\infty}M_{\rho_n}
=\underline{u}_\gamma(\bar{x})-\overline{u}_\gamma(\bar{x}),$
which contradicts \eqref{061}. Consequently, $B_2\geq D_2$ does  not hold.\\

(ii) Now, we look at the case $B_1\geq D_1$. Then we  have
\begin{eqnarray}
e^{\gamma x_{\rho_n}}(\gamma \underline{u}_\gamma(x_{\rho_n})+\zeta_{{\rho_n}}^\prime (x_{\rho_n}))\leq
e^{\gamma y_{\rho_n}}(\gamma \overline{u}_\gamma(y_{\rho_n})+\varphi_{{\rho_n}}^\prime (y_{\rho_n})).\label{0a1}
\end{eqnarray}

In the rest of the proof we consider $x_{\rho_n}\ge b_0$ and $y_{\rho_n}\ge b_0$ only.
It follows immediately from  \eqref{a12} and \eqref{0a1} that
\begin{eqnarray}
&&e^{\gamma x_{\rho_n}} \underline{u}_\gamma(x_{\rho_n})-e^{\gamma y_{\rho_n}}\overline{u}_\gamma(y_{\rho_n})
 \leq \frac{\frac{4m}{\left({\rho_n}(y_{\rho_n}-x_{\rho_n})^2+1\right)^2}-\rho_n}{\gamma}(y_{\rho_n}-x_{\rho_n}) (e^{\gamma y_{\rho_n}}-e^{\gamma x_{\rho_n}})
 .\label{a8}
\end{eqnarray}
Let $N_2(\epsilon)$ be a positive integer such that for all $n\geq N_2(\epsilon)$, $\rho_n\geq 4m$.  Since $(y_{\rho_n}-x_{\rho_n})(e^{ry_{\rho_n}}-e^{rx_{\rho_n}})$ is always nonnegative, then from \eqref{a8} we can see  that for all $n\geq N_2(\epsilon)$,
\begin{eqnarray}
e^{\gamma x_{\rho_n}} \underline{u}_\gamma(x_{\rho_n})-e^{\gamma y_{\rho_n}}\overline{u}_\gamma(y_{\rho_n})\leq 0.\label{0a4}
\end{eqnarray}

Recall that $x_{\rho_n}\rightarrow \overline{x}$, $y_{\rho_n}\rightarrow \overline{y}$ and $\overline{x}= \overline{y}$. There exists an integer $N_3(\epsilon)$ such that for all $n\geq N_3(\epsilon)$,
\begin{eqnarray*}
|e^{\gamma x_{\rho_n}}-e^{\gamma \overline{x}}|<\epsilon, \ \ |e^{\gamma y_{\rho_n}}-e^{\gamma \overline{x}}|<\epsilon \ \ \mbox{and}\ \ |\overline{u}_\gamma( x_{\rho_n})-\overline{u}_\gamma( y_{\rho_n})|<\epsilon.
\end{eqnarray*}
Then for $n\geq N_3(\epsilon)$, we have
\begin{eqnarray}
&&\underline{u}_\gamma(x_{\rho_n})(1-e^{\gamma x_{\rho_n}})-\overline{u}_\gamma(y_{\rho_n})(1-e^{\gamma y_{\rho_n}})\nonumber\\
&=&\underline{u}_\gamma(x_{\rho_n})(1-e^{\gamma x_{\rho_n}})-\overline{u}_\gamma(x_{\rho_n})(1-e^{\gamma y_{\rho_n}})+(\overline{u}_\gamma(x_{\rho_n})-\overline{u}_\gamma(y_{\rho_n}))(1-e^{\gamma y_{\rho_n}})\nonumber\\
&<& \underline{u}_\gamma(x_{\rho_n})(1-e^{\gamma \overline{x}}+\epsilon)-\overline{u}_\gamma(x_{\rho_n})(1-e^{\gamma \overline{x}}-\epsilon)+\epsilon\nonumber\\
&\leq&M
(1-e^{\gamma\overline{x}})+(\underline{u}_\gamma(x_{\rho_n})
+\overline{u}_\gamma(x_{\rho_n})+1)
\epsilon,\label{0a3}
\end{eqnarray}
where the lat inequality follows by \eqref{cc15}.

\noindent
Since the functions $\underline{u}_\gamma$ and $\overline{u}_\gamma$ are bounded, it can be easily shown that \begin{eqnarray}
\frac{Me^{\gamma \overline{x}}}{\sup_{x}|\underline{u}_\gamma(x)
+\overline{u}_\gamma(x)+1|}>0.\label{7a9}
\end{eqnarray}
 From \eqref{3080}, \eqref{3084}, \eqref{0a4}, \eqref{0a3} and \eqref{7a9}, it follows  that
 for any $\epsilon<\frac{Me^{\gamma \overline{x}}}{\sup_{x}|\underline{u}_\gamma(x)+\overline{u}_\gamma(x)+1|}$ and $n\geq \max\{N_2(\epsilon),N_3(\epsilon)\}$,
\begin{eqnarray*}
M\leq M_{\rho_n}+\frac{2m}{\rho_n}&=& \underline{u}_\gamma(x_{\rho_n})-\overline{u}_\gamma(y_{\rho_n})\\
&=& e^{\gamma x_{\rho_n}}\underline{u}_\gamma(x_{\rho_n})-e^{\gamma y_{\rho_n}}\overline{u}_\gamma(y_{\rho_n})\\
&&+\underline{u}_\gamma(x_{\rho_n})(1-e^{\gamma x_{\rho_n}})-\overline{u}_\gamma(y_{\rho_n})(1-e^{\gamma y_{\rho_n}})\\
&<&M
(1-e^{\gamma\overline{x}})+(\underline{u}_\gamma(x_{\rho_n})
+\overline{u}_\gamma(x_{\rho_n})+1)
\epsilon
< M,
\end{eqnarray*}
which is an contradiction. So $B_1\geq D_1$ does not hold.\\

Combining (i) and (ii) shows that $B_1<D_1$  and $B_2<D_2$. This is a contraction to the fact that  at least of the inequalities $B_1\geq D_1$  and $B_2\geq D_2$ holds. As a result, $\underline{u}(x)\leq \overline{u}(x)$ for all $x\geq -\frac p\alpha$.
This completes the proof.
\hfill $\square$\\

\bigskip

\noindent {\bf Proof of Theorem \ref{T3}} Assume that $x\in (-\frac p\alpha,\overline{x}]$.

i) Let $\Pi(n)$ denote the set of admissible strategies such that if the initial reserve $x<\bar{x}$, the controlled reserve will always stay below or at $\bar{x}$ until the arrival of the $n$th claim.

We will show that for any $n\in \mathbb{N}$ and $x\in(-\frac p\alpha,\overline{x}]$, $V(x)=\sup_{L\in \Pi(n)}V_L(x)$ by induction.

Noting that $\Pi(0)=\Pi$, we get $V(x)=\sup_{L\in \Pi(0)}V_L(x).$

Assume that $V(x)=\sup_{L\in \Pi(n-1)}V_L(x)$ for some $n\geq 1$.

Let $L^{(n-1,x)}\in \Pi(n-1)$ be an $\frac \epsilon2$-optimal strategy for the reserve process with  the initial value $x$, that is
\begin{eqnarray}
0\leq V(x)-V_{L^{(n-1,x)}}(x)\leq \frac \epsilon2.\label{00b5}
\end{eqnarray}

Let  $\tau^L $ denote the first time that the reserve process under strategy $L$ reaches $\bar{x}$, and $\hat{\tau}^L$  the arrival time of the next claim occurring  after time $\tau^L$.

Then given  the initial reserve $x$ ($x\leq\bar{x}$),  we can construct an $\frac \epsilon2$-optimal strategy  $L^{(n,x)}\in \Pi({n})$ as follows. Apply the strategy $L^{(n-1,x)}$ until the first time the controlled reserve reaches $\bar{x}$, then pays out dividends at a rate equal to the premium incoming rate to keep the reserve at the level $\bar{x}$ until the arrival of the next claim. After that, we apply the strategy $L^{(n-1,R_{\hat{\tau}^{L(n,x)}})}$ to the shifted process $\theta_{\hat{\tau}^{L^{(n,x)}}}R$.

Recall that $S_1$ and $U_1$ are respectively the arrival time and the amount of the first claim. Note that for the case with initial reserve $R_0=\bar{x}$, under strategy $L^{(n,x)}$ we have
$$\hat{\tau}^{L^{(n,x)}}=S_1,\ \ R_{S_1}=(\overline{x}-U_1)\vee (-\frac p\alpha) \ \ \mbox{ and }$$
$$R_t=\bar{x}, \ \  \ \dif L^{(n,\overline{x})}_t= (p+r\bar{x}I\{\bar{x}\geq 0\}+\alpha \bar{x}I\{\bar{x}<0\})\dif t \mbox{ for $t<S_1$}.
$$
Hence,  by noticing the fact that ruin will not occur before the arrival of the first claim, i.e. $T\geq S_1$, and that $V(-\frac p\alpha)=0$, we obtain that  given the initial reserve $\bar{x}$,
\begin{eqnarray}
V_{L^{(n,\bar{x})}}(\bar{x})&=&\E_{\bar{x}}\Big[ \int_0^{S_1}e^{-\delta s}\left(p+r\bar{x}I\{\bar{x}\geq 0\}+\alpha \bar{x}I\{\bar{x}<0\}\right)\dif s\nonumber\\
&&+e^{-\delta S_1}V_{L^{\left(n-1,(\bar{x}-U_1)\vee (-\frac p\alpha)\right)}}\left((\bar{x}-U_1)\vee (-\frac p\alpha)\right)\Big]\nonumber\\
&=&\frac 1{\lambda+\delta}\left(p+r\bar{x}I\{\bar{x}\geq 0\}+\alpha \bar{x}I\{\bar{x}<0\}\right)\nonumber\\
&&+\frac \lambda{\lambda+\delta}\int_0^{\bar{x}+\frac p\alpha}V_{L^{(n-1,\bar{x}-y)}}(\bar{x}-y)\dif F(y)\label{7a20}
\end{eqnarray}
It follows by \eqref{7a21}, \eqref{00b5}, \eqref{7a20} and  assumption that $\mathcal{G}_V(\overline{x})=0$ that
\begin{eqnarray}
V_{L^{(n,\bar{x})}}(\bar{x})&\geq&\frac 1{\lambda+\delta}\left(p+r\bar{x}I\{\bar{x}\geq 0\}+\alpha \bar{x}I\{\bar{x}<0\}\right)\nonumber\\
&&+\frac \lambda{\lambda+\delta}\int_0^{\bar{x}+\frac p\alpha}(V(\bar{x}-y)-\frac \epsilon2)\dif F(y)\nonumber\\
&\geq & V(\bar{x})-\frac \epsilon2.\label{b12}
\end{eqnarray}
  Note that for any fixed $x\in[-\frac p\alpha,\bar{x}]$ and for $k=n-1$ and $n$, we have
\begin{eqnarray}
V_{L^{(k,{x})}}({x})
&=&\E_x\left[ \int_0^{\tau^{L^{(k,{x})}}}e^{-\delta s}\dif  L^{(k,{x})}_s;\tau^{L^{(k,{x})}}<T\right]\nonumber\\
&&+\E_x[ e^{-\delta \tau^{L^{(k,{x})}}};\tau^{L^{(k,{x})}}<T]
V_{L^{(k,\bar{x})}}(\bar{x})+\E_x[ e^{-\delta \tau^{L^{(k,{x})}}};\tau^{L^{(k,{x})}}\geq T].\label{7a22}
\end{eqnarray}
From the construction of the strategies, we can see that given the initial reserve $x$,
\begin{eqnarray}
\tau^{L^{(n,{x})}}=\tau^{L^{(n-1,{x})}}\ \ \ \mbox{and}\ \ \ L_s^{(n,x)}=L_s^{(n-1,x)} \ \mbox{for } s\leq \tau^{L^{(n,{x})}}.\label{7a23}\end{eqnarray}
By using \eqref{7a22} for $k=n-1$ and $n$, and \eqref{7a23}, we obtain
\begin{eqnarray}
V_{L^{(n-1,{x})}}({x})-V_{L^{(n-1,{x})}}({x})=\E_x\left[ e^{-\delta \tau^{L^{(n,{x})}}};\tau^{L^{(n,{x})}}<T\right]
\left(V_{L^{(n,\bar{x})}}
(\bar{x})-V_{L^{(n-1,\bar{x})}}(\bar{x})\right).\nonumber\\
\label{b11}
\end{eqnarray}
Note that by the definition of $L^{(n-1,{x})}$ and \eqref{b12} we have
 \begin{eqnarray}
 V(\bar{x})\geq V_{L^{(n-1,\bar{x})}}(\bar{x})\geq V(\bar{x})-\frac \epsilon2\ \ \mbox{and } \ V(\bar{x})\geq V_{L^{(n,\bar{x})}}(\bar{x})\geq V(\bar{x})-\frac \epsilon2,\label{b13}
 \end{eqnarray}
which implies
\begin{eqnarray}
V_{L^{(n,\bar{x})}}(\bar{x})- V_{L^{(n-1,\bar{x})}}(\bar{x})\geq -\frac \epsilon2.\label{00b6}
\end{eqnarray}
Combining \eqref{00b5}, \eqref{b11} and \eqref{00b6} gives
 $V_{L^{(n,{x})}}({x})\geq V({x}) -\epsilon.$
Therefore, $$V(x)\geq \sup_{L\in\Pi(n)}V_L(x)\geq V_{L^{(n,{x})}}({x})\geq V({x}) -\epsilon.$$
Consequently, letting $\epsilon\rightarrow 0$ gives us $V({x})=\sup_{L\in\Pi(n)}V_L(x) \ \mbox{ for $x\in [-\frac p\alpha,\overline{x}]$.}$

ii) Now we try to find a strategy $\hat{L}\in\Pi_{\overline{x}}$ such that it is $\epsilon$-optimal.

 Noting that $V(\overline{x})>0$, we can find a $t_1$  large enough such that
\begin{eqnarray}
e^{-\delta t_1}<\frac{\epsilon}{4V(\bar{x})}.\label{cc20}
\end{eqnarray}
Then for this fixed $t_1$, choose an $n$ large enough such that
\begin{eqnarray}\pr(N(t_1)\geq n)=\sum_{k\geq n}\frac{e^{-\lambda t_1}(\lambda t_1)^k}{k!}\leq \frac{\epsilon}{4V(\bar{x})}.\label{cc21}
\end{eqnarray}

Define $\sigma^L$ to be the first time that the controlled reserve process under strategy $L$ reaches $\overline{x}$ after the arrival of the $n$th claim ($S_n$).

Let $L^{(n,x)}$ be any $\frac \epsilon2$-optimal strategy in $\Pi(n)$ given the initial reserve $x$.  Given the initial reserve $x$, construct an dividend payout strategy $\hat{L}(x)$ such that the strategy  $L^{(n,x)}$ is applied before time $\sigma^{L^{(n,x)}}$, then at time $t=\sigma^{L^{(n,x)}}$, a lump sum of $\bar{x}+\frac p\alpha$ is paid out immediately,   and thereafter no dividends will be paid out.

Then, we have \begin{eqnarray}
V_{\hat{L}({x})}({x})
&=&\E_x\left[ \int_0^{\sigma^{L^{(n,x)}}}e^{-\delta s}\dif  L^{(n,{x})}_s;\sigma^{L^{(n,x)}}<T\right]+\E_x[ e^{-\delta \sigma^{L^{(n,x)}}};\sigma^{L^{(n,x)}}<T](\bar{x}+\frac p\alpha)\nonumber\\
&&+\E_x\left[ \int_0^{T}e^{-\delta s}\dif  L^{(n,{x})}_s;\sigma^{L^{(n,x)}}\geq T\right].\label{b14}
\end{eqnarray}
 Note that for any initial reserve $x\leq\bar{x}$, the strategy $\hat{L}(x)$ is same as $L^{(n,{x})}$  until both the controlled reserve under the former strategy reaches $\bar{x}$ for the first time, which implies $$\sigma^{\hat{L}(x)}=\sigma^{L^{(n,x)}} \ \ \mbox{ and }\ \     \hat{L}_{t}(x)=L^{(n,x)}_t \ \ \mbox{ for $t\leq \sigma^{L^{(n,x)}}$}.  $$

  Noting that \eqref{7a22} also holds for the strategy $L^{(n,x)}$ here, by \eqref{b14} we get, for $x\leq\bar{x}$,
\begin{eqnarray}
V_{\hat{L}({x})}({x})-V_{L^{(n,{x})}}({x})&=&
\E_x[ e^{-\delta \sigma^{L^{(n,x)}}};\sigma^{L^{(n,x)}}<T](\bar{x}+\frac p\alpha-V_{L^{(n,\overline{x})}}(\bar{x}))\nonumber\\
&\geq& -\E_x[ e^{-\delta \sigma^{L^{(n,x)}}}]V(\bar{x})
.\label{0b18}
\end{eqnarray}

  As $S_n\leq \sigma^{L^{(n,x)}}$, we have $\{\sigma^{L^{(n,{x})}}<t_1\}\subseteqq \{S_n<t_1\} \subset \{N(t_1)\geq n\}$ for $x\leq\bar{x}$. Therefore, for any $x\leq\bar{x}$,
\begin{eqnarray}
\{\sigma^{L^{(n,x)}}<\infty\} \subset \{\sigma^{L^{(n,{x})}}\geq t_1\} \cup \{N(t_1)\geq n\}.\label{0b20}
\end{eqnarray}
Then by \eqref{cc20} and \eqref{cc21} we have
\begin{eqnarray}
\E_x[ e^{-\delta \sigma^{L^{(n,x)}}}]&\leq& \E_x[ e^{-\delta \sigma^{L^{(n,x)}}};\sigma^{L^{(n,{x})}}\geq t_1]+\E_x[ e^{-\delta \sigma^{L^{(n,x)}}};N(t_1)\geq n]\nonumber\\
&\leq& e^{-\delta t_1}+\pr(N(t_1)\geq n)\leq \frac{\epsilon}{2V(\bar{x})}.\label{0b19}
\end{eqnarray}
It follows from \eqref{0b18} and \eqref{0b19} that for $x\in[-\frac p\alpha,\overline{x}]$
\begin{eqnarray}
V_{\hat{L}({x})}({x})\geq V_{L^{(n,{x})}}({x})-\frac\epsilon2\geq V(x)-\epsilon,
\end{eqnarray}
where the last inequality is due to the fact that $L^{(n,x)}$ is an $\frac \epsilon2$-optimal strategy.

Noting that $\hat{L}(x)\in\Pi_{\overline{x}}$, the above inequality implies  \begin{eqnarray*}
\sup_{L\in \Pi_{\bar{x}}}V_{L}({x})\geq V(x),\ \ \ x\in[-\frac p\alpha, \overline{x}].
\end{eqnarray*}
This concludes the proof.
\hfill $\square$\\

\bigskip

\noindent {\bf Proof  of Theorem \ref{T4}} It is sufficient to show that for any $\epsilon>0$, there exists a strategy $\overline{L}^{(x)}\in \Pi_{\bar{x}}$ such that
\begin{eqnarray}
V_{\overline{L}^{(x)}}({x})\geq V(x)-\epsilon \mbox{ for all $x\in (-\frac p\alpha,\bar{x}]$.}\label{b22}
\end{eqnarray}

For a positive $\epsilon<4V(\bar{x})$, define
\begin{eqnarray}
\Delta(\epsilon)=\frac{p+r\overline{x}I\{\overline{x}\geq 0\}+\alpha \overline{x}I\{\overline{x}< 0\}}{\delta}\ln\frac{4V(\bar{x})}{\epsilon},\label{cc43}
\end{eqnarray}
\begin{eqnarray}
x_n&=&\overline{x}-\frac{\Delta(\epsilon)}{n}, \ \mbox{and }\ h_n=\frac{V(x_n)-V(\bar{x})}{x_n-\bar{x}}-1. \label{cc40}
\end{eqnarray}
It can be shown that $x_n\leq \overline{x}$ and $x_n\rightarrow \overline{x}$.
Since $V^\prime(\bar{x})=1$, we have  $\lim_{n\rightarrow \infty}h_n=0$. Moreover, notice
\begin{eqnarray*}
&&\lim_{n\rightarrow\infty}\left(\frac{rx_{n}+p}{r\overline{x}+p}\right)^{\frac{\delta n}{r}}I\{\overline{x}\geq  0\}+\left(\frac{\alpha x_{n}+p}{\alpha\overline{x}+p}\right)^{\frac{\delta n}{\alpha}}I\{\overline{x}< 0\}\\
&=& \exp\left\{\frac{-\delta \Delta(\epsilon)}{p+r\overline{x}I\{\overline{x}\geq 0\}+\alpha \overline{x}I\{\overline{x}< 0\}}\right\}.
 \end{eqnarray*}
 Hence we can choose a $n_0$ such that
\begin{eqnarray}
&&\left(\frac{rx_{n_0}+p}{r\overline{x}+p}\right)^{\frac{\delta n_0}{r}}I\{\overline{x}\geq 0\}+\left(\frac{\alpha x_{n_0}+p}{\alpha\overline{x}+p}\right)^{\frac{\delta n_0}{\alpha}}I\{\overline{x}<0\}\nonumber\\
&\leq& \exp\left\{\frac{-\delta \Delta(\epsilon)}{p+r\overline{x}I\{\overline{x}\geq 0\}+\alpha \overline{x}I\{\overline{x}< 0\}}\right\}+\frac{\epsilon}{4V(\overline{x})},\label{cc46}
\end{eqnarray}
 and \begin{eqnarray}
h_{n_0}<\frac{\epsilon}{8\Delta(\epsilon)}.\label{7a25}
\end{eqnarray}

For any $x\geq -\frac p\alpha$, let $L^{(0,x)}$ be a $\frac{\epsilon}{8n_0}$-optimal strategy given the initial reserve $x$, that is
\begin{eqnarray}
V_{L^{(0,x)}}({x})\geq V(x)-\frac{\epsilon}{8n_0} . \label{b25}
\end{eqnarray}

Let  $\tau^L $ denote the first time that the controlled reserve process under strategy $L$ reaches $\bar{x}$ starting from an initial reserve below $\bar{x}$.

For $x\leq \overline{x}$, define a sequence of strategies $\{L^{(n,x)}\}_{n\geq 1}$ recursively as follows:  $L^{(n,x)}$
is a strategy given the initial reserve $x$ that the insurer pays dividends according to strategy $L^{(n-1,x)}$ until the reserve reaches $\overline{x}$ for the first time ($\tau^{L^{(n-1,x)}}$),  pays out a lump sum of $\overline{x}-x_{n_0}$ at time $\tau^{L^{(n-1,x)}}$, and thereafter employs the strategy $L^{(n-1,x_{n_0})}$ to the shifted process $ \theta_{\tau^{L^{(n-1,x)}}}R$.

 It can be shown that for all $n$, $\tau^{L^{(n,x)}}=\tau^{L^{(0,x)}}$, and $L^{(n,x)}_s=L^{(0,x)}_s$ for $s\leq \tau^{L^{(0,x)}}$. Then we have for $x\in(-\frac p\alpha,\overline{x}]$ and $n=1,2,\cdots$,
\begin{eqnarray}
V_{L^{(n,{x})}}({x})
&=&\E_x\left[ \int_0^{\tau^{L^{(0,x)}}}e^{-\delta s}\dif  L^{(0,{x})}_s;\tau^{L^{(0,x)}}<T\right]\nonumber\\
&&+\E_x[e^{-\delta \tau^{L^{(0,x)}}};\tau^{L^{(0,x)}}<T]\left(V_{L^{(n-1,{x_{n_0}})}}({x_{n_0})}
+\overline{x}-x_{n_0}\right)\nonumber\\
&&+\E_x\left[ \int_0^T e^{-\delta s}\dif  L^{(n-1,{x})}_s;\tau^{L^{(0,{x})}}>T\right].\label{b23}
\end{eqnarray}
Using  \eqref{b25} for $x=x_{n_0}$ and $\bar{x}$,  \eqref{b23} for $n=1$ and the second equality in \eqref{cc40}, we get
\begin{eqnarray}
&&|V_{L^{(1,x)}}({x})-V_{L^{(0,x)}}(x)|\nonumber\\
&=&| \E_x[e^{-\delta \tau^{L^{(0,{x})}}};\tau^{L^{(0,{x})}}<T]
\left(\overline{x}-x_{n_0}+V_{L^{(0,x_{n_0})}}(x_{n_0})
-V_{L^{(0,\overline{x})}}(\overline{x})\right)|\nonumber\\
&\leq&|\overline{x}-x_{n_0}-V(\overline{x})
+V(x_{n_0})|+V(x_{n_0})-V_{L^{(0,x_{n_0})}}(x_{n_0})+V(\overline{x})
-V_{L^{(0,\overline{x})}}(\overline{x})\nonumber\\
&\leq &h_{n_0}(\overline{x}-x_{n_0})+\frac{\epsilon}{4n_0}\leq \frac{3\epsilon}{8n_0},\label{b26}
\end{eqnarray}
where the last inequality follows by the first equality in \eqref{cc40} and \eqref{7a25}.\\
Therefore, from \eqref{b23} we have for $x\in(-\frac p\alpha,\overline{x}]$ and $n\geq 2$,
\begin{eqnarray}
|V_{L^{(n,{x})}}({x})-V_{L^{(n-1,{x})}}({x})|&\leq& \E_x[e^{-\delta \tau^{L^{(0,{x})}}}]|V_{L^{(n-1,x_{n_0})}}(x_{n_0})
-V_{L^{(n-2,{x_{n_0}})}}({x_{n_0}})|\nonumber\\
&\leq&|V_{L^{(1,x_{n_0})}}(x_{n_0})
-V_{L^{(0,{x_{n_0}})}}({x_{n_0}})|\leq\frac{3\epsilon}{8n_0}.\label{b26}
\end{eqnarray}
Consequently, by \eqref{b25} and \eqref{b26}
\begin{eqnarray*}
|V(x)-V_{L^{(n_0,{x})}}({x})|&=&|V(x)-V_{L^{(0,{x})}}({x})
+\sum_{n=1}^{n_0}\left(V_{L^{(n-1,{x})}}({x})-V_{L^{(n,{x})}}({x})\right)|
\\&<& \frac\epsilon{8n_0}+\frac{3\epsilon}{8n_0}\leq \frac{\epsilon}{2}.
\end{eqnarray*}

Define $\bar{\tau}=\inf\{t>0:R_t^{L^{(n_0,x)}}>\bar{x}\}$, where $R_t^{L^{(n_0,{x})}}$ represent the controlled reserve process under strategy $L^{(n_0,{x})}$.

Under strategy $L^{(n_0,{x})}$, in order to exceed $\bar{x}$, the controlled reserve process with initial reserve $x_{n_0}$ should go from $x_{n_0}$ up to $\bar{x}$ for at least $n_0$ times. Note from the dynamics \eqref{dyn} that it will take at least $t_0(x_{n_0},\bar{x})$ (defined in \eqref{cc2})
  for this reserve process to reach $\bar{x}$ starting from $x_{n_0}$.
Therefore, $\bar{\tau}\geq n_0t_0(x_{n_0},\bar{x}) $. Consequently,  it follows by \eqref{cc2} \eqref{cc43}, \eqref{cc40} and \eqref{cc46} that
\begin{eqnarray}
\E_{x_{n_0}}[e^{-\delta \bar{\tau}}]\leq \E[e^{-\delta  n_0t_0(x_{n_0},\bar{x})}] \leq \frac\epsilon{2V(\bar{x})}.\label{b28}
\end{eqnarray}

Next, we construct a strategy $\bar{L}(x)$ through $L^{(n_0,{x})}$: pays dividends according to the strategy $L^{(n_0,{x})}$ before time $\bar{\tau}$ (the time that the reserve process reaches $\overline{x}$ for the first time), pays out a lump sum of $\bar{x}+\frac p\alpha$ at time  $\bar{\tau}$, and thereafter pays no dividends.

\noindent Then we have
\begin{eqnarray}
V_{\bar{L}({x})}({x})&=&\E_x[ \int_0^{\bar{\tau}} e^{-\delta s}\dif  \bar{L}_s({x});\bar{\tau}<T]+\E_x[ e^{-\delta \bar{\tau}};\bar{\tau}<T](\bar{x}+\frac p\alpha)\nonumber\\
&&+\E_x[ \int_0^{T} e^{-\delta s}\dif  \bar{L}_s({x});\bar{\tau}>T].\label{b29}
\end{eqnarray}
Notice that
\begin{eqnarray}
V_{L^{(n_0,{x})}}({x})&=&\E_x[ \int_0^{\bar{\tau}} e^{-\delta s}\dif  L^{(n_0,x)}_s;\bar{\tau}<T]+\E_x[ \int_{\bar{\tau}}^T e^{-\delta s}\dif  L^{(n_0,{x})}_s;\bar{\tau}<T]\nonumber\\
&&+\E_x[ \int_0^T e^{-\delta s}\dif  L^{(n_0,x)};\bar{\tau}>T],
\label{b30}
\end{eqnarray}
and
\begin{eqnarray}
\E_x[ \int_{\bar{\tau}}^T e^{-\delta s}\dif  L^{(n_0,{x})}_s;\bar{\tau}<T]\leq \E[e^{-\delta \bar{\tau}}]V(\bar{x}).\label{b31}
\end{eqnarray}
Since  $\bar{L}_s(x)=L^{(n_0,x)}_s$ for $s\leq \bar{\tau}$, it follows from  \eqref{b29}, \eqref{b30} and \eqref{b31} that for  $x\in(-\frac p\alpha,\overline{x}]$,
\begin{eqnarray*}
V_{\bar{L}(x)}(x)-V_{L^{(n_0,{x})}}({x})\geq \E_x[e^{-\delta \bar{\tau}}]\left(\bar{x}+\frac p\alpha-V(\bar{x})\right)\geq -\E_x[e^{-\delta \bar{\tau}}]V(\bar{x})\geq -\frac\epsilon2,
\end{eqnarray*}
where the last inequality is due to \eqref{b28}.

So $\bar{L}(x)$ is the desired strategy.
\hfill $\square$\\

\bigskip

\noindent {\bf Proof of Lemma \ref{L0c}}

(a) Since $\Lambda_{V}(x )$ is continuous in $x$, $\mathcal{A} $ is closed.\\

(b) (i) To prove that $\mathcal{B} $ is left-open, it is sufficient to show that for any $x\in \mathcal{B} $ we can find an $h>0$ such that for any $ y\in(x-h,x) $, $V^\prime(y )<1$.

Note that $V^\prime(x )=1$ for $x\in \mathcal{B} $ and
$G_{x-h }^\prime(y )=1$, and  that $p+ryI\{y\geq 0\}+\alpha yI\{y<0\}$ is increasing in $y$. Therefore, it follows from \eqref{3} that for any $y\in(x-h,x)$,
\begin{eqnarray}
\mathcal{L}_{G_{x-h }}(y )&\leq&\mathcal{L}_V(x )
-(\lambda+\delta)(V(x)-G_{x-h }(y))
+\nonumber\\
&&\lambda \left(\int_{0}^{y+\frac p\alpha}G_{x-h}(y-u )\dif F (u)-\int_{0}^{x+\frac p\alpha}V(x-u )\dif F (u)\right). \label{JEaa7}
\end{eqnarray}
    Noticing $\mathcal{L}_V(x )<0$, $G_{x-h }(y )=V(y )$ for $y\leq x-h$, and $\lim_{h\rightarrow 0}G_{x-h }(y )\rightarrow V(x )$ for $y\in(x-h,x)$, it follows from \eqref{JEaa7} that
    \begin{eqnarray}
\mathcal{L}_{G_{x-h }}(y )&<& 0  \mbox{ for small $h>0$}.  \label{JEaa8}
\end{eqnarray} This along with the fact that $G_{x-h}'(y )\equiv 1$ for $y>x-h$ implies that $G_{x-h }$ is a viscosity super-solution to \eqref{1}  on $(x-h,x]$. Then by Theorem \ref{theorem39} (ii), we have $V(y )=G_{x-h }(y )$ for all $y \in [-\frac p\alpha,x]$. As a result,
\begin{eqnarray}
V'(y ) = G_{x-h }'(y )=1, \mbox{ for } y\in (x-h,x]. \label{JEaa9}
\end{eqnarray}
Combining \eqref{JEaa8} and \eqref{JEaa9} implies $(x-h,x]\subset \mathcal{B} .$
Therefore, $\mathcal{B} $  is left-open.

(ii) To prove that there exist a $y$ such that $(y,\infty)\subset \mathcal{B} $, it is sufficient to show that we can find a large enough $y>0$ such that $\mathcal{L}_{G_{y }}(x )< 0$ for all $x>y$, because  if $G_{y }(x )$ of this kind is a super-solution on $(y,\infty)$ and therefore $V'(x ) = G_{y }'(x )\equiv 1$ for $x>y$.

Noticing that $G_{y }(x )$ is nondecreasing in $x$, and that  $G_{y }(x )=x-y+V(y )$ and  $V(y )-y>\frac p\alpha$ for $x>y$, we obtain that for $x>y>0$,
\begin{eqnarray*}
\mathcal{L}_{G_{y}}(x )&=&p+rx-(\lambda +\delta)(G_{y }(x )+\lambda \int_{0}^{x+\frac p\alpha}G_{y}(x-y )\dif F (y)\\
&<&p-(\delta-r ) x-\delta y -\delta\frac p\alpha<0 \mbox{ for large $y$}
,\label{0aa96}
\end{eqnarray*}
where the last inequality follow by noticing $\delta>r$.

The existence of $y$ also indicates that $\mathcal{B} $ is not empty.

\vskip 0.75cm

\noindent (c) Noticing that $V(-\frac p\alpha)=0$, by the definition of $\mathcal{G}$ \eqref{7a21}, we obtain $\mathcal{G}_V(-\frac p\alpha )=0$, which implies $-\frac p\alpha\in \mathcal{A} $.

Assume that $x_1>x_0>-\frac p\alpha$, $(x_0,x_1]\subset \mathcal{B} $ and $x_0 \notin  \mathcal{B} $. We will show in the following that $x_0\in \mathcal{A} $.

If $V'(x_0 )=1$, then from the fact that $x_0 \notin  \mathcal{B} $ and $\mathcal{L}_V\leq 0$, we know $\mathcal{G}_V(x_0 )=0$. Therefore, $x_0\in \mathcal{A} $

Now assume, on the other hand,  $V'(x_0 )\neq 1$. It follows from the fact $(x_0,x_1]\subset \mathcal{B} $ that $V'(x )=1$ for all $x\in (x_0,x_1]$, which  implies
\begin{eqnarray}
\lim_{x\downarrow x_0}\frac{V(x )-V(x_0 )}{x-x_0}=1. \label{0aa81}
\end{eqnarray}

 Define $$a=\liminf_{x\uparrow x_0}\frac{V(x )-V(x_0 )}{x-x_0}.$$
   By Lemma \ref{thm2.2}  we know $a\geq 1$. We distinguish two cases: 1. $a>1$ and 2. $a=1$.

Case 1: Assume $a>1$. Then for any $b$ with $1<b\leq a$, we have
$$\limsup_{x\downarrow x_0}\frac{V(x )-V(x_0 )}{x-x_0}=1<b<\liminf_{x\uparrow x_0}\frac{V(x )-V(x_0 )}{x-x_0}.$$
Since $V$ is a viscosity sub-solution, by Remark \ref{cc50} (i), it follows that there exists a continuously differentiable function $\phi: (-\frac p\alpha,\infty)\rightarrow \mathbb{R}$ such that $V-\phi$ reaches a maximum at $x_0 $ with $\phi^\prime(x_0 )=b$. Therefore, by Definition \ref{alt} (i) it follows that
\begin{eqnarray*}
&&\max\left\{1-b,(p+rx_0I\{x_0\geq 0\}+\alpha x_0 I\{x_0<0\}) b-(\lambda +\delta)V(x_0 )\right.\\
&&+\left.\lambda \int_{0}^{x_0+\frac p\alpha}V(x_0-y )\dif F (y)\right\}\geq 0,
\end{eqnarray*}
which implies
 $$(p+rx_0I\{x_0\geq 0\}+\alpha x_0I\{x_0<0\}) b-(\lambda+\delta)V(x_0 )+\lambda \int_{0}^{x_0+\frac p\alpha}V(x_0-y )\dif F (y)\geq 0.$$
Taking limits $b\rightarrow 1$ gives $\mathcal{G}_V(x_0 )\geq 0.$
Since $\mathcal{G}_V(x )$ is continuous in $x$ and $\mathcal{G}_V(x )<0$ for $x\in(x_0,x_1]$, it can be seen that $\mathcal{G}_V(x_0 )= 0$, which implies $x_0\in \mathcal{A} $.\\

Case 2: Assume $a=1$.  then we can find a sequence $\{h_n\}$ with $h_n\downarrow 0$ such that
\begin{eqnarray}
\lim_{n\rightarrow \infty} \frac{V(x_0 )-V(x_0-h_n )}{h_n}=1. \label{JEaa12}
\end{eqnarray}
Define $$a_n=\frac{V(x_0 )-V(x_0-h_n )}{h_n}-1,\ \mbox{and}\ A_n=\{x\in [0,h_n]: V'(x ) \mbox{ exists and } V'(x )\geq 1+2a_n\}.$$
By Theorem \ref{thm2.2} we know that $a_n\geq 0.$

i) If there exists some $n$ such that $a_n=0$, then we have
\begin{eqnarray}
V(x_0 )-V(x )=x_0-x \ \mbox{ for $x\in[x_0-h_n,x_0]$}. \label{JEaa10}
\end{eqnarray}
Otherwise, if for some $x\in[x_0-h_n,x_0]$, $V(x_0 )-V(x )>x_0-x$, then
\begin{eqnarray*}
V(x_0 )-V(x_0-h_n )&=&V(x_0 )-V(x )+V(x )-V(x_0-h_n )\\
 &>&x_0-x+x-(x_0-h_n)=h_n,
 \end{eqnarray*}
which contradicts the assumption $a_n=0$.

As a result of \eqref{0aa81} and \eqref{JEaa10}, we have $V'(x_0 )=1.$ Therefore,  $\mathcal{G}_V(x_0 )\geq 0$ follows by noticing $x_0 \in \mathcal{B} $. Notice that $\mathcal{G}_V(x_0 )\leq 0$ due to the continuity of $\mathcal{G}_V$. Therefore, $\mathcal{G}_V(x_0 )=0$, implying $x_0 \in \mathcal{A} $.

ii) Suppose  $a_n>0$ for all $n$. Since $V(x )$ is differentiable almost everywhere, and $V^\prime(x )$, if exists, is greater than $1$, we have
\begin{eqnarray*}
a_n+1&=& \frac{\int_{0}^{h_n}V'(x )\dif x}{h_n}=\frac{\int_{ A_n}V'(x )\dif x+\int_{[0,h_n]\setminus A_n} V'(x ) \dif x}{h_n}\\
&\geq& \frac{|A_n|(1+2a_n)+(h_n-|A_n|)}{h_n}\label{ja1},
\end{eqnarray*}
where $|A_n|$ denotes the Lebesgue measure of the set $A_n$. It follows from \eqref{ja1} that $|A_n|\leq \frac{h_n}{2}\rightarrow 0$.
Therefore we can find a sequence $x_n\uparrow x_0$ such that $V'(x_n )$ exist and $1\leq V'(x_n )< 1+2a_n$.
Consequently, $\lim_{n\rightarrow \infty}V'(x_n )=1.$

If there exists a subsequence $\{x_{nj}\}$ with $x_{nj}\uparrow x_0$ such that $V'(x_{nj} )>1$, then by \eqref{1} we have $\mathcal{G}_V(x_{nj} )=0$. This implies $x_{nj}\in \mathcal{A} $. Since $\mathcal{A} $ is a closed set, we conclude that $x_0\in \mathcal{A} $.

Suppose that there is an integer $n_0>0$, such that for all $n\geq n_0$, $V'(x_n )=1$. We will show by Proof by Contradiction that $\mathcal{G}_V(x_0 )=0.$
 Assume $\mathcal{G}_V(x_0 )<0$.
Let $n$ be large enough such that
\begin{eqnarray}
V(x_0 )-V(x_n )<-\mathcal{G}_V(x_0 )/(\lambda+\delta). \label{JEaa13}
\end{eqnarray}
Note that $V(y)\geq V(x_n )+y-x_n=G_{x_n}(x )$ for all $y\geq x_n$. Then for all $x\in [x_n,x_0]$,
\begin{eqnarray}
\mathcal{G}_{G_{x_n}}(x )&=&p+rxI\{x\geq 0\}+\alpha x I\{x<0\}-(\lambda+\delta)G_{x_n}(x )
+\lambda \int_{0}^{x+\frac p\alpha}G_{x_n}(x-y )\dif F (y)\nonumber\\
&\leq& \mathcal{G}_V(x_0 )+(\lambda+\delta)(V(x_0 )-(V(x_n )+x_0-x_n))\nonumber\\
&\leq& \mathcal{G}_V(x_0 )+(\lambda +\delta)(V(x_0 )-V(x_n ))< 0,
\label{JEaa14}
\end{eqnarray}
where the last inequality follows from \eqref{JEaa13}.

\noindent Noting that $G_{x_n}'(x )=1$ for $x>x_n$ and $G_{x_n}(x )=V_{x_n}(x )$ for $x\in [0,x_n]$, so by \eqref{JEaa14} it follows
$\mathcal{L}_{G_{x_n}}(x )=\mathcal{G}_{G_{x_n}}(x )< 0 \mbox{ for all } x\in [x_n,x_0].$ Therefore $G_{x_n}$ is a viscosity super-solution on $[x_n,x_0]$.

\noindent Recalling that $\mathcal{G}_V(x_0)=0$, then by Theorem \ref{theorem39} (ii), we have
$V(x )=G_{x_n}(x ) \mbox{ for } x\in [0,x_0].$
As a result, $V(x )$  is differentiable at $x_0$ and
\begin{eqnarray}
 V'(x_0 )=G_{x_n}'(x_0 )=1. \label{JEaa15}
\end{eqnarray}
Combining \eqref{JEaa14} and \eqref{JEaa15}  implies $x_0 \in \mathcal{B} $, which contradicts the fact that $x_0\notin \mathcal{B} .$ Therefore $\mathcal{G}_{V}(x_0 )\geq 0$.\\
Since $V$ is a viscosity super-solution and $V'(x_0 )=1$, from the definition of viscosity super-solution we can see that
$\mathcal{G}_{V}(x_0 )=\mathcal{L}_{V}(x_0 )\leq 0. $\\
Consequently, $\mathcal{G}_V(x_0 )=0$, which implies $x_0\in \mathcal{A} $.\\

\noindent (d) For any $x\in \mathcal{C} $, we have $\mathcal{G}_V(x )<0$. Since $\mathcal{G}_V(x )$ is continuous, we can find a $\epsilon$ small enough such that
\begin{eqnarray}
\mathcal{G}_V(y )<0 \mbox{ for all $y\in [x,x+\epsilon)$}.\label{new2}
\end{eqnarray}

If for all $y\in(x,x+\epsilon)$, $y\notin \mathcal{B} $, then $[x,x+\epsilon)\subset \mathcal{C}  $.

If, on the other hand, there exist an $x_1\in (x,x+\epsilon)$ such that $x_1\in \mathcal{B} $, then we can find an $x_0$ and $x_1$ with $x_0<x_1$ such that $x_0\in \mathcal{A} $ and  $(x_0,x_1]\subset \mathcal{B} $.
As $x<x_1$ and $x\notin \mathcal{B} $, we conclude that $x_0\in (x,x_1)\subset (x,x+\epsilon)$,
 which along with $\mathcal{G}_V(x_0 )=0$ is a contradiction to \eqref{new2}.
This completes the proof.
\hfill $\square$\\


\bigskip
\bigskip

\end{document}